\documentclass[12pt, a4paper, parskip=half, abstracton]{scrartcl}
\pdfoutput=1

\usepackage{array}
\usepackage{marginnote}
\usepackage{xcolor}
\usepackage{amscd,amssymb,amsfonts,amsmath,latexsym,amsthm}
\usepackage{hyperref}
\usepackage[all,cmtip]{xy}
\usepackage{mathrsfs}
\usepackage{graphicx}
\usepackage{bm} 
\usepackage{etoolbox}
\def\hB{\hspace*{\fill}$\qed$}
\usepackage{relsize} 

\usepackage{mathtools} 

\usepackage[nottoc]{tocbibind} 
\usepackage{defs_pp1}
\usepackage{slashed}
\usepackage[utf8]{inputenc}
\usepackage{microtype}
\usepackage[english]{babel}

\usepackage[titles]{tocloft} 
\title{The coarse index class with support}

\author{
Ulrich Bunke\thanks{Fakult{\"a}t f{\"u}r Mathematik,
Universit{\"a}t Regensburg,
93040 Regensburg,
GERMANY\newline
\href{mailto:ulrich.bunke@mathematik.uni-regensburg.de}{ulrich.bunke@mathematik.uni-regensburg.de}} 
\and
Alexander Engel\thanks{Institut f{\"u}r Mathematik und Informatik der Universit{\"a}t Greifswald\newline
Walther-Rathenau-Str.\,47 in 17489 Greifswald, GERMANY\newline
\href{mailto:alexander.engel@uni-greifswald.de}{alexander.engel@uni-greifswald.de}}
}
\numberwithin{equation}{section}
\newtheorem{theorem}{Theorem}[section] 
\newtheorem{prop}[theorem]{Proposition}
\newtheorem{lem}[theorem]{Lemma}
\newtheorem{ddd}[theorem]{Definition}
\newtheorem{kor}[theorem]{Corollary}

\theoremstyle{remark}
\theoremstyle{definition}

\newtheorem{ex}[theorem]{Example}
\newtheorem{rem}[theorem]{Remark}

\newcommand{\Res}{\mathrm{Res}}

\newcommand{\Hilb}{\mathbf{Hilb}}

\newcommand{\inter}{\mathrm{int}}

\newcommand{\BC}{\mathbf{BornCoarse}}

\newcommand{\Fib}{{\mathrm{Fib}}}

\newcommand{\incl}{\mathrm{incl}}

\newcommand{\bF}{{\mathbf{F}}}

\newcommand{\bG}{{\mathbb{G}}}

\newcommand{\cY}{{\mathcal{Y}}}

\renewcommand{\epsilon}{\varepsilon}

\newcommand{\KX}{K\!\mathcal{X}}

\renewcommand{\Dirac}{\slashed{D}}

\newcommand{\dist}{\mathrm{dist}}

\newcommand{\IR}{\mathbb{R}}

\newcommand{\Ccat}{{\mathbf{C}^{\ast}\mathbf{Cat}}}
\newcommand{\nCcat}{{\mathbf{C}^{\ast}\mathbf{Cat}}^{\mathrm{nu}}}
\newcommand{\grnCcat}{{}^{\mathrm{gr}}{\mathbf{C}^{\ast}\mathbf{Cat}}^{\mathrm{nu}}}

\newcommand{\nCalg}{{\mathbf{C}^{\ast}\mathbf{Alg}}^{\mathrm{nu}}}
\newcommand{\grnCalg}{{}^{\mathrm{gr}}{\mathbf{C}^{\ast}\mathbf{Alg}}^{\mathrm{nu}}}

\newcommand{\on}{\mathit{on}}
\newcommand{\Spin}{\mathit{Spin}}
\newcommand{\spin}{\mathit{spin}}
\newcommand{\lc}{\mathit{lc}}

\begin{document}

\maketitle

\vspace*{-2ex}
\begin{abstract}
\noindent
We construct the coarse index class  with a support condition   of an equivariant Dirac operator on a complete Riemannian manifold endowed with a proper and isometric action of a group as   a class of the spectrum-valued equivariant coarse $K$-homology theory $K\!\cX^{G}$ from \cite{bu}. Moreover we show a coarse relative index theorem and discuss the compatibility of the index  with the suspension isomorphism. 
\end{abstract}

\tableofcontents

\section{Introduction}
 
%

Coarse (co)homology theories were introduced by Roe \cite{roe_index_1,MR1147350} in order to produce large-scale invariants of  metric  spaces. Since then {they} experienced a broad development with many applications to index theory, geometric group theory,  isomorphism conjectures like the Baum--Connes and Farrell--Jones conjectures, surgery theory or the classification of metrics of positive scalar curvature.
 
In order to provide a very general  formal framework  for coarse geometry capturing all known examples in  \cite{equicoarse} 
we introduced   the category $G\BC$ of $G$-bornological coarse spaces. These are $G$-sets equipped with a compatible  bornology and coarse structure. Metric spaces with isometric $G$-action represent $G$-bornological coarse spaces in a natural way. 
We further proposed a notion of an equivariant coarse homology theory and studied the universal example of such a homology theory.
In the present paper we in particular  consider the equivariant coarse $K$-homology $\KX^{G}$ which is constructed in detail in \cite{bu}. 

The applications of the  category $G\BC$ and coarse homology theories to the Farrel-Jones and Novikov conjectures have been developed in \cite{unik}, \cite{desc} and \cite{fajo}. The goal of the present paper make the general theory and the functor $K\cX^{G}$ applicable to index theory and to provide the background for  \cite{coarsevew}, \cite{Bunke:2021aa}.

Indeed, index theory  was one of the original sources of interest in coarse geometry. 
Roe \cite{MR1147350} observed that one can define an index of a Dirac operator on a complete Riemannian manifold as a $K$-theory class of a certain $C^{*}$-algebra, now called the Roe algebra, capturing large-scale properties of the  underlying Riemannian manifold. This index is interesting since it encodes non-trivial information on the non-invertibility of the Dirac operator.
The main objective of the present paper  is to provide a reference for the construction of   the index class of a $G$-equivariant Dirac operator on a complete Riemannian manifold $M$ with isometric $G$-action as a class in  $K\cX^{G}(M)$, i.e., the value   of the coarse homology theory $K\cX^{G}$ on the bornological coarse space represented by $M$. This is non-trivial since, in contrast to classical definitions, $K\cX^{G}(M)$ is not defined as the $K$-theory of the specific Roe algebra carrying the index. This will be explained in this introduction in greater detail further below.


One possible reason for the invertibility of a Dirac operator and hence the triviality of the index is the positivity of the zero-order term in the Weizenboeck formula,  e.g., the scalar curvature of the manifold in the case of the spin Dirac operator. If this term is   positive   on some subset  of the manifold only, then the  index information is supported on the complement of that subset. The precise formulation of this fact is due to Roe \cite{Roe:2012fk} and leads to the notion of the coarse index class of a Dirac operator with support.

 Let $M$ be a complete Riemannian manifold with a proper action of a discrete group $G$ by isometries and let $A$ be an invariant subset of $M$. If $\Dirac $ is a  $G$-equivariant (degree $n$)-Dirac operator   which is uniformly locally positive outside of $A$, then 
    we construct an index class \begin{equation}\label{iojofijwefewofefew}
\Ind(\Dirac ,\on(A))\:\:\mbox{in}\:\: \KX^{G}_{n}(A)\ .
\end{equation} For this
  {we need the technical assumption that $A$ is a support. It requires essentially that a cofinal subset of the invariant coarse thickeninings of $A$ are coarsely equivalent to $A$. For a precise formulation we refer to  Definition \ref{rgklherogervfvvfgrege}}.

 In order to deal with the degree $n$ in a simple and uniform manner, following Zeidler \cite{MR3551834}, we use operators commuting with an action of the Clifford algebra $\Cl^{n}$. We   follow this reference further  by representing the index class by a homomorphism from the graded algebra $\hat C_{0}(\R)$ to the Roe algebra obtained by functional calculus applied to the Dirac operator.

   In the present paper we  show the following basic properties of this index class:
  \begin{enumerate}
  \item  locality: Coarse Relative Index Theorem \ref{roigeorgergergerg}
  \item  compatibility with suspension: Theorem \ref{fiweofwefwefewf}
\end{enumerate}

%
%
%
%
%
In \cite{coarsevew} we show that  this index class together with the knowledge of  the above two basic  properties  are sufficent to apply the  theory of  abstract boundary value problems to Dirac operators. This allows to  reproduce the main results of Piazza--Schick  \cite{MR3286895} and Zeidler \cite{MR3551834} 
on secondary coarse index  invariants by calculations in   coarse homotopy theory without 
opening the black-box of analysis of Dirac operators anymore.

 Let $\Dirac $ be an equivariant Dirac operator on $M$ which  is uniformly locally positive outside of an invariant subset $A$ of $M$.   
The first construction of an  index class with support in the non-equivariant case {was given by Roe} \cite{Roe:2012fk}. In this reference the index is an element of the $K$-theory of a Roe algebra which in Roe's notation is written as $C^{*}(A\subseteq X)$. 

For manifolds with a free and proper $G$-action such an index class has been constructed by Piazza--Schick \cite{MR3286895} and Zeidler \cite[Def.~4.5]{MR3551834}. Zeidler's index class is an element of the $K$-theory of Yu's localization algebra denoted  in Zeidler's notation by  $C^{*}_{L,A}(M)^{G}$. His index class is actually a refinement of the index classes considered in the present paper and by Piazza--Schick. The latter  correspond to the image of Zeidler's index class under the natural evaluation homomorphism $C^{*}_{L,A}(M)^{G}\to C^{*}(A\subseteq X)^{G}$. The superscript  indicates the  subalgebra of  operators in the Roe algebra $C^{*}(A\subseteq M)$  which can be approximated by $G$-invariant finite propagation operators. 

 The assumption made by  Zeidler that the action of the group is free leads to the simplification that the controlled Hilbert space obtained from the $L^{2}$-sections of the bundle on which the Dirac operator acts is already ample. For proper actions this is not the case and a stabilization is necessary. In the present paper we in particular  also provide the generalization of the construction of the  index class from free to proper actions.
Further applications are given in \cite[Sec. 6]{Bunke:2021aa}

The main goal of the present paper is to explain how the index class of a Dirac operator of degree $n$ on $M$
can be captured as a class in $K\cX^{G}_{n}(M)$. In order to explain why this is not obvious we recall the classical construction of  equivariant coarse $K$-homology  and compare it with the construction of $K\cX^{G}$.

Classically the equivariant coarse $K$-homology is a $\Z$-graded group valued functor 
defined on suitably nice metric spaces with isometric $G$ which  sends coarse equivalences to isomorphisms and  has a Mayer-Vietoris sequence for equivariant coarsely excisive decompositions. 
 In   \cite{equicoarse}  took the properties of such group-valued functors as the basis for an axiomatization of  the notion of an equivariant coarse homology theory  as a spectrum-valued functor defined on the category of $G$-bornological coarse spaces. The functor $K\cX^{G}$ from \cite{bu} is an example of such an equivariant coarse homology theory in the sense of  \cite{equicoarse} which by  Theorem \ref{fwefiwjfeiooi234234324434e}    indeed refines the 
  classical group-valued functor.
 
 We first recall the
 classical    construction, see e.g.  \cite{willett_yu_book}.  Given a sufficiently nice metric space with isometric $G$-action $X$ one choses a representation of $X$ consisting of a separable Hilbert space $H$ with unitary $G$-action  and  a covariant representation $\phi :C_{0}(X)\to B(H)$
which is sufficiently large (ample). This choice gives rise to the Roe algebra $C(X,H,\phi)$ consisting of the
  bounded operators on $H$ which can be approximated by invariant  locally compact finite propagation operators. One then 
 defines the coarse $K$-homology of $X$ as the $K$-theory groups  $K_{*}( C(X,H,\phi))$ of the Roe algebra.
One  checks that these groups are independent on the choice of $(H,\phi)$. If $f:X\to X'$ is an equivariant  controlled and proper map and $(H',\phi')$ is a   choice of a representation  for $X'$, then one can choose
 a controlled equivariant  isometric inclusion $U:(H ,f_{*}\phi)\to (H',\phi')$ which induces a homomorphism
 $U-U^{*}:C(X,H,\phi)\to C(X',H',\phi')$ of Roe algebras. One then checks that the induced map of $K$-groups does not depend on the choices. 
 
 In the case of  a Dirac operator $\Dirac$ on a bundle $E\to M$ for $M$ a complete Riemannian manifold of non-zero dimension with free action of $G$   one can take $H:=L^{2}(M,E)$ and $\phi:C_{0}(M)\to B(H)$ the action by multiplication operators. Then the index of $\Dirac$ is naturally a class in the $K$-theory of the Roe algebra  $C(M,H,\phi)$ of the ample $(H,\phi)$,
  or even the subalgebra $C(A\subseteq M,H,\phi)$ if $\Dirac $ is uniformly locally positive outside $A$.
 
 We now put this in contrast to the construction of $K\cX^{G}$. In order to ensure spectrum-valued functoriality we 
avoid the necessity of choosing ample representations which anyway do not exist in the generality of arbitrary $G$-bornological coarse spaces. Instead we functorially associate to  any $G$-bornological coarse  space $X$ the $C^{*}$-category  $\bC(X)$ of equivariant locally finite $X$-controlled Hilbert spaces  and then define $K\cX^{G}(X):=K(\bC(X))$.
The details of this construction have been worked out in \cite{bu}  and will be recalled in Section~\ref{kfopkwefewdqwedcacsdcac}.
 
 In view of this discussion our main non-trivial task (achieved in Section \ref{goijrgoiregregreg}) is to construct for a complete Riemannian manifold with proper $G$-action and equivariant bundle
 a canonical map
 $$K_{*}(C(M,H,\phi))\to K_{*}(\bC(X))\ ,$$ where $H:=L^{2}(M,E)$ and $\phi:C_{0}(M)\to B(H)$ is as above. 
 This map will then be used to transport the index class of $\Dirac$ in $K_{*}(C(M,H,\phi))$ to a class in 
 $K\cX^{G}_{*}(X)=K_{*}(\bC(X))$.
 In order to incorporate the degree appropriately, in Section \ref{erkogpqwerfrewfrefrefw} we develop a version of the above constructions
 for graded Roe algebras and graded $C^{*}$-categories.

\begin{ex}
Here is the typical example. Let $M$ be an $n$-dimensional Riemannian spin manifold with a proper action of $G$ by isometries  and assume that the scalar curvature is bounded below by a positive constant on the subset $M\setminus A$ where $A$ is a support. {Denote by} $P\to M$ the $G$-equivariant $\Spin(n)$-principal bundle. We can consider $\Spin(n)$ as a subgroup of the units in $\Cl^{n}$ which acts on $\Cl^{n}$ by left-multiplication.
Then we define the Dirac bundle $E:=P\times_{\Spin(n)}\Cl^{n}$ with the $\Cl^{n}$-action given by right-multiplication. 
Let $\Dirac_{M}^{\spin}$ be the  associated $G$-invariant $\Cl^{n}$-linear   Dirac operator on $E$. Then we get a class
\[\Ind(\Dirac^{\spin}_{M}, {\on(A)})\:\:\mbox{in}\:\: \KX^{G}_{n}(A )\ .\]
 \hB
\end{ex}

 
%
%
%


\paragraph{Acknowledgements}
\textit{The authors were supported by the SFB 1085 ``Higher Invariants'' funded by the Deutsche Forschungsgemeinschaft DFG.}

\textit{The second named author was further supported by a Research Fellowship of the Deutsche Forschungsgemeinschaft DFG (EN 1163/1-1 ``Mapping Analysis to Homology'', project number 314131583) which he spent at the Texas A\&M University.}

\textit{The first named author thanks Thomas Schick for motivating discussions.}

\section{\texorpdfstring{$\boldsymbol{X}$}{X}-controlled Hilbert spaces} 
 
In this section we introduce the notion  of an equivariant    $X$-controlled Hilbert space and the notions of controlled and of locally compact operators. To keep the paper self-contained we reproduce some definitions from \cite{buen} and \cite{bu}.

  Let $G$ be a discrete group. A $G$-Hilbert space is a Hilbert space with a unitary action  of $G$. We will usually use a letter like $H$ to denote a $G$-Hilbert space and   the notation  $(g,h)\mapsto gh$ for the action. A $G$-$*$-algebra is a $*$-algebra  with an action of $G$ by automorphisms of $*$-algebras. 
 If the underlying algebra is a $C^{*}$-algebra, then we say that it is a $G$-$C^{*}$-algebra.
  For example,  the algebra of bounded operators $B(H)$ on a $G$-Hilbert space $H$ is a $G$-$C^{*}$-algebra with the action by conjugation 

If $A$ is a $G$-$*$-algebra, then the subset of $G$-invariants $A^{G}$ is a $*$-algebra. Similarly, if $A$ is a $G$-$C^{*}$-algebra, then $A^{G}$ is a $C^{*}$-algebra.  
 
%

 We  consider a $G$-bornological coarse   space $X$. It is a $G$-set  with a $G$-coarse structure $\cC$  and
 a bornology $\cB$. The elements of $\cC$ are  subsets of $X\times X$ and called the coarse  entourages of $X$, while the  elements of $\cB$ are subsets of $X$ and called the bounded subsets of $X$. 
  \begin{ddd}\label{riferofergrgreg}
 An equivariant $X$-controlled Hilbert space is a pair $(H,\phi)$ of a Hilbert space with a unitary action of $G$ and a finitely-additive projection-valued measure $\phi$ defined on all subsets of $X$ such that $H(B)$ is separable for every bounded subset of $X$ and  $\phi$ is equivariant.
 \end{ddd}
Here we use the notation $H(Y):=\phi(Y)(H)$ for the subspaces of $H$ corresponding to subsets $Y$ of $X$. 
Equivariance of $\phi$ means that
 for every $g$ in $G$ and subset $Y$ of $X$ we have the equality of projections
$$g\phi(Y) g^{-1}=\phi(g^{-1}Y)\ .$$
Note that the notion of an equivariant $X$-controlled Hilbert space only depends on the underlying $G$-bornological space of $X$ and not on the coarse structure.


Let
  $(H,\phi)$ and $(H^{\prime},\phi^{\prime})$ be two equivariant $X$-controlled Hilbert spaces and $A:H\to H^{\prime}$ be a bounded operator.
\begin{ddd}\mbox{}
\begin{enumerate}
\item $A$ is locally compact if for every bounded subset $B$ of $X$ the operators
$\phi^{\prime}(B)A$ and $A\phi(B)$ are compact.
\item $A$ is controlled if there exists an entourage $U$ of $X$ such that
$\phi^{\prime}(Y^{\prime})A\phi(Y)=0$ for every  two subsets  $Y$ and $Y^{\prime}$  of $X$ such that $U[Y]\cap Y^{\prime}=\emptyset$.
\item $A$ is invariant, if $gAg^{-1}=A$ for every $g$ in $G$.
\end{enumerate}
\end{ddd}

Assume now that $X$ is a $G$-bornological coarse space.
Let furthermore $(H,\phi)$ and $(H^{\prime},\phi^{\prime})$ be two equivariant $X$-controlled Hilbert spaces.

\begin{ddd}\label{fiuoffewfewfweff}
We let $C_{lc}(X,(H^{\prime},\phi^{\prime}),(H,\phi))$ denote  closure (in operator norm) of the set of   invariant bounded   operators  from $H^{\prime}$ to $H$
 which are  locally compact and controlled.
\end{ddd}



%


%

The collection of spaces 
$C_{lc}(X,(H^{\prime},\phi^{\prime}),(H,\phi))$ 
 for all pairs of equivariant   $X$-controlled Hilbert spaces is closed under 
 taking adjoints or  forming linear combinations or  compositions.
 In particular
$$C_{lc}(X,H,\phi):=C_{lc}(X,(H,\phi),(H,\phi))
$$
is a $C^{*}$-algebra.

\begin{rem}\label{fweoihoewfweffwe}
Note that   $C_{lc}(X,(H^{\prime},\phi^{\prime}),(H,\phi))$ consists of $G$-invariant operators. We will refrain from using the notation $C_{lc}(X,(H^{\prime},\phi^{\prime}),(H,\phi))^{G}$ since this is misleading. We would use this notation for the $G$-invariants in the Banach space  $C_{lc}(\Res^{G}_{1}(X,(H^{\prime},\phi^{\prime}),(H,\phi)))$ (this means that we first forget the $G$-action on the data, then form the Banach space of operators in the non-equivariant case, and finally take the invariants under the $G$-action). The space $C_{lc}(X,(H^{\prime},\phi^{\prime}),(H,\phi))$ might be strictly contained in $C_{lc}(X,(H^{\prime},\phi^{\prime}),(H,\phi))^{G}$. 
\hB
\end{rem}

\section{Locally finite \texorpdfstring{$\boldsymbol{X}$}{X}-controlled Hilbert spaces}

In order to define equivariant coarse $K$-homology  $\KX^{G}(X)$
we consider the $C^{*}$-category of equivariant $X$-controlled Hilbert spaces with the morphism spaces defined in Definition~\ref{fiuoffewfewfweff}. In order to ensure that this $C^{*}$-category is unital we must restrict to 
locally finite   equivariant $X$-controlled Hilbert spaces. In this section we recall this notion and prove some technical results which will be crucial later.

Let $X$ be a $G$-bornological coarse space.

\begin{ddd} An equivariant $X$-controlled Hilbert space is determined on points if  
the natural map $ \bigoplus_{x\in X} H(\{x\})\to H$ is an isomorphism.
\end{ddd}
The sum above is understood in the sense of Hilbert spaces, i.e., it involves a completion.

Let $(H,\phi)$ be an equivariant $X$-controlled Hilbert space.

\begin{ddd} The support $\supp(H,\phi)$ of  $(H,\phi)$ is the  subset
$\{x\in X\:|\: H(\{x\})\not=0\}$ of the space $X$.
\end{ddd}

Note that the support of an equivariant $X$-controlled Hilbert space
  is an invariant subset of $X$.

 {Recall that a subset $S$ of $X$ is called locally finite, if the intersection $S \cap B$ is finite for every bounded subset $B$ of $X$.}

\begin{ddd}
An equivariant $X$-controlled Hilbert space is  locally finite if the following conditions are satisfied:
\begin{enumerate}
\item $(H,\phi)$ is determined on points.
\item The support  of $(H,\phi)$   is locally finite.  
\item For every $x$ in $X$ the space $H(\{x\})$ is finite-dimensional.
\end{enumerate}
\end{ddd}

 If $(H,\phi)$ is locally finite, then for every bounded subset $B$ of $X$ the space $H(B)$ is finite-dimensional.

 Let $(H,\phi)$ be an equivariant  $X$-controlled Hilbert space, not necessarily locally finite.
\begin{ddd}
A closed invariant subspace $\tilde H$ of $H$ is called locally finite if there exists a locally finite equivariant 
$X$-controlled Hilbert space $(\tilde H^{\prime},\phi^{\prime})$ and a controlled equivariant  isometric embedding
$U: {\tilde H^{\prime}}\to H$ with image $\tilde H$.
\end{ddd}

Let $A$ be a bounded operator on $H$. 
\begin{ddd}
We say that $A$ is locally finite, if there exist locally finite subspaces
$\tilde H_{0}$ and $\tilde H_{1}$ of $H$ and an operator $\tilde A:\tilde H_{0}\to \tilde H_{1}$ such that
$A=j_{1}\circ \tilde A\circ p_{0}$, where $p_{0}$ is the  {orthogonal} projection from $H$ to $\tilde H_{0}$ and
$j_{1}$ is the inclusion of $\tilde H_{1}$ into $H$. 
\end{ddd}

Note that a locally finite operator is automatically locally compact.

\begin{ddd}
We let $C(X,H,\phi)$  be the closure of the subset of  locally finite operators in $C_{lc}(X,H,\phi)$.
\end{ddd}

 With the structures induced from $C_{lc}(X,H,\phi)$ the subset $C(X,H,\phi)$  is a $C^{*}$-algebra.

We now consider the possibility of approximating locally compact operators by  locally finite ones. Our approach uses the existence of sufficiently nice equivariant partitions of the space $X$.

Let $X$ be a $G$-set.
An equivariant family of subsets $(Y_{i})_{i\in I}$ of $X$ is a family of subsets where $I$ is a $G$-set and we have the equality $gY_{i}=Y_{gi}$ for every $i$ in $I$ and $g$ in $G$.  

Let $X$ be a $G$-bornological coarse space.
 Recall that  the action of $G$ on $X$ is called proper if for every bounded subset $B$ the set $\{g\in G\:|\: gB\cap  B\not=0\}$ is finite. Properness only involves the bornology of $X$. 
For our comparison results we need a stronger condition which we call \emph{very proper}. 

Let $X$ be a $G$-bornological coarse space.

\begin{ddd}\label{fpowjoefewfewf}
We say that $X$ is very proper if for every entourage $U$ of $X$ there exists an  equivariant   partition $(B_{i})_{i\in I}$ of $X$   such that
\begin{enumerate}
\item\label{iegergerg} The partition is uniformly bounded,  {i.e., there exists an entourage $V$ of $X$ such that every subset $B_i$ is $V$-bounded}, i.e., $B_{i}\times B_{i}\subseteq V$,
\item \label{hweiofjweiofewf} $I$   has finite stabilizers,
\item  \label{hweiofjweiofewf1} for every $i$ in $I$ the set $\{j\in I\:|\: U[B_{i}]\cap B_{j}\not=\emptyset\}$ is finite,
  \item\label{ojreoigjoergergerg} for every $i$ in $I$ there exists a point  in $B_{i}$ which is fixed by the stabilizer of $i$,
  \item\label{gergihieorgergrege} for every bounded subset $B$ of $X$ the set
$\{i\in I\:|\: B_{i}\cap B\not=\emptyset\}$ is finite, and
\item \label{hweiofjweiofewf5} the set of orbits $I/G$ is countable.
\end{enumerate}
 {If $X$ also has the structure of a $G$-measurable space,} we say that $X$ is measurably very proper if
we in addition can assume that the members of the partition are measurable.
\end{ddd}

 {The following proposition ensures that under suitable conditions the usual definition of the Roe algebra using locally compact operators coincides with the definition that we use which employs locally finite operators.}

Let $X$ be a $G$-bornological coarse space and $(H,\phi)$ an equivariant  $X$-controlled Hilbert space.

\begin{prop}\label{lem:sdfbi23}
If $X$ is  very proper and $(H,\phi)$ is determined on points, then the natural inclusion is an equality $C(X,H,\phi) = C_{\lc}(X,H,\phi)$.
\end{prop}

\begin{proof} 
We consider an operator  $A$ in $C_{\lc}(X,H,\phi)$. We will show that for every given $\epsilon$ in $(0,\infty)$ there exists a  locally finite  
subspace $H^{\prime}$ of $H$  
and an operator $A^{\prime}:H^{\prime}\to H^{\prime}$ in $C_{\lc}(X,H,\phi)$ such that $\|A- A^{\prime}  \|\le \epsilon$.
Here we omit to write the projection from $H$ to $H^{\prime}$  in front of $A^{\prime}$ and the inclusion of $H^{\prime}$ into $H$ after $A^{\prime}$.

As a first step we choose an invariant  locally compact operator $A_{1}$ of propagation controlled by an entourage $\supp(A_{1})$ such that $\|A-A_{1}\|\le \epsilon/2$. This is possible by the definition of $C_{\lc}(X,H,\phi)$.

%
%

Since $X$ is very proper   we can choose an   equivariant partition $(B_{i})_{i\in I}$ of $X$
with the properties listed in Definition \ref{fpowjoefewfewf}  for the invariant entourage $\supp(A_{1})$ in place of $U$. 

We assume that $I/G$ is infinite (the argument in the finite case is similar, but simpler).
Then by Assumption \ref{fpowjoefewfewf}.\ref{hweiofjweiofewf5} we can assume that $I=\bigsqcup_{n\in \nat} I_{n}$ for transitive $G$-sets $I_{n}$.
For every integer $n$ we choose a base point $i_{n}$ in $I_{n}$.

For every integer $n$ let $P_{i_{n}}$ be a finite-dimensional $G_{i_{n}}$-invariant projection on $H(B_{i_{n}})$, where $G_{i_{n}}$ is the subgroup of $G$ fixing $i_{n}$. We will fix this projection later.  Since $B_{i_{n}}$ is $G_{i_{n}}$-invariant we can approximate the identity of $H(B_{i_{n}})$ strongly by such projections. Here we use the Assumption \ref{fpowjoefewfewf}.\ref{hweiofjweiofewf} which implies that $G_{i_{n}}$ is finite.
For $i$ in $I_{n}$ we define $P_{i}:=gP_{i_{n}}g^{-1}$, where $g$ in $G$ is such that $gi_{n}=i$. This projection on $H(B_{i})$ is well-defined, and the family of projections $(P_{i})_{i\in I}$ is $G$-invariant in the sense that $gP_{i}g^{-1}=P_{gi}$ for every $i$ in $I$ and $g$ in $G$.

We now define the $G$-invariant operator
$$A^{\prime}:=\sum_{i,j\in I} P_{i}A_{1}P_{j}\ .$$
Then
\begin{align*}
A_{1}-A^{\prime} & = \sum_{i,j\in I}\phi(B_{i})A_{1}\phi(B_{j})-\sum_{i,j\in I}P_{i}A_{1}P_{j}\\
& = \sum_{i,j\in I}(\phi(B_{i})-P_{i})A_{1}\phi(B_{j})+\sum_{i,j\in I}P_{i}A_{1}(\phi(B_{j})-P_{j})\ .
\end{align*}
For every $j$ in $I$  we define the set $$J_{j}:= \{i\in I\:|\: B_{i}\cap \supp(A_{1})[B_{j}]\not=\emptyset\}\ .$$
This set is finite by Assumption \ref{fpowjoefewfewf}.\ref{hweiofjweiofewf1}.
Then we have
$$A_{1}-{A^\prime}=\sum_{j}\sum_{i\in J_{j}}(\phi(B_{i})-P_{i})A_{1}\phi(B_{j})+\sum_{j}\sum_{i\in J_{j}}P_{i}A_{1}(\phi(B_{j})-P_{j})\ .$$
Using the orthogonality of the terms for different indices $i$
we get the estimates
$$\|A_{1}-{A^\prime}\| \le \max_{i\in I} \sum_{\{j\in I|i\in J_{j}\}}\|(\phi(B_{i})-P_{i})A_{1}\|+\max_{i\in I}
  \sum_{\{j\in I| i\in J_{j}\}} \|A_{1}(\phi(B_{j})-P_{j})\|\ .$$
The index sets of these sums are finite by Assumption \ref{fpowjoefewfewf}.\ref{hweiofjweiofewf1}. For $i$ in $I$ we set
$$\ell_{i}:=|\{j\in I|i\in J_{j}\}|\ .$$
 
 Since $A_{1}$ is locally compact we can make $\|A_{1}(\phi(B_{i})-P_{i})\|$ and $\|(\phi(B_{i})-P_{i})A_{1}\|$ as small as we want by choosing $P_{i}$ sufficiently big.
 
  If we choose for every integer $n$ the projection
 $P_{i_{n}}$ so large that
 \begin{equation}\label{flwehfjoiewfefewf}\|(\phi(B_{i_{n}})-P_{i_{n}})A_{1}\|\le \frac{\epsilon}{4\ell_{i_{n}}}\ ,\end{equation}
 then we get (using $G$-invariance)
$$ \max_{i\in I} \sum_{\{j\in I|i\in J_{j}\}}\|(\phi(B_{i})-P_{i})A_{1}\|\le \frac{\epsilon}{4}\ .$$
In order to deal with the second term we must increase the projections further. We proceed by induction.
Assume that we have choosen the projections such that
$$\max_{i\in G\{i_{0},i_{1},\dots,i_{n}\}}  \sum_{\{j\in I| i\in J_{j}\}} \|A_{1}(\phi(B_{j})-P_{j})\|\le
 \frac{\epsilon}{4}\ .$$
In the next step we further increase the projections
$P_{i_{m}}$ for all integers $m$  with $  i_{n+1}\in G J_{i_{m}}$ (these are finitely many by Assumption \ref{fpowjoefewfewf}.\ref{hweiofjweiofewf1})
such that
$$\max_{i\in G\{i_{0},i_{1},\dots,i_{n+1}\}}  \sum_{\{j\in I| i\in J_{j}\}} \|A_{1}(\phi(B_{j})-P_{j})\|\le
 \frac{\epsilon}{4}$$
 and \eqref{flwehfjoiewfefewf} is satisfied.
We now observe that with this procedure we increase every projection at most a finite number of times (by Assumption \ref{fpowjoefewfewf}.\ref{hweiofjweiofewf1}).

Let $H^{\prime}$ be the sub-Hilbert space of $H$ spanned by the images of $P_{i}$ for all $i$ in $I$. Then $A^{\prime}$ factorizes over $H^{\prime}$ and 
$$\|A-A^{\prime}\|\le \epsilon\ .$$ By  \ref{fpowjoefewfewf}.\ref{iegergerg} can choose an entourage $V$ which bounds  the members of the  partition.
Then the propagation of $A^{\prime}$ is bounded by $V\circ \supp(A_{1})\circ V$.



We now construct a control $\phi^{\prime}$ for $H^{\prime}$ exhibiting this subspace as a locally finite subspace.
For every integer $n$ we choose a point $b_{i_{n}}$ in $B_{i_{n}}$ which is fixed by $G_{i_{n}}$. This is possible by Assumption \ref{fpowjoefewfewf}.\ref{ojreoigjoergergerg}. For $i$ in $I_{n}$ we
then define $b_{i}:=gb_{i_{n}}$ where $g$ in $G$ is such that $gi_{n}=i$. Note that $b_{i}$ is well-defined. The collection of points $(b_{i})_{i\in I}$ thus defined is equivariant in the sense that
$gb_{i}=b_{gi}$ for every $i$ in $I$ and $g$ in $G$.

We now define the equivariant projection-valued measure
$$\phi^{\prime}:=\sum_{i\in I} \delta_{b_{i}}P_{i}$$ on $H^{\prime}$.
It turns $(H^{\prime},\phi^{\prime})$
into an equivariant $X$-controlled Hilbert space which is determined on points.
It is locally finite   by Assumption \ref{fpowjoefewfewf}.\ref{gergihieorgergrege}  and the fact that $P_{i}$ is finite-dimensional for every $i$ in $I$.

Finally we see that the inclusion
$H^{\prime}\to H$ is $V$-controlled.
Hence $H^{\prime}$ is a locally finite subspace of $H$.
\end{proof}

Let $X$ be a $G$-bornological coarse space
and $(H,\phi)$ be an equivariant $X$-controlled Hilbert space.
Let $H^{\prime}$ and $H^{\prime\prime}$ be two   subspaces of $H$.
Then we define
$$H^{\prime}\bar{+} H^{\prime\prime}:=\overline{H^{\prime}+H^{\prime\prime}}\ .$$
{The following rather technical result will be needed in the proof of Theorem~\ref{fwefiwjfeiooi234234324434e}.}
\begin{prop}\label{iofeoewfewfewf}
Assume that
\begin{enumerate}
\item $X$ is very proper,
\item $(H,\phi)$ is determined on points, and
\item $H^{\prime}$ and $H^{\prime\prime}$ are locally finite.
\end{enumerate}
Then $H^{\prime}\bar{+}H^{\prime\prime}$ is locally finite.
\end{prop}
\begin{proof}
By assumption we can find a symmetric invariant entourage $U$ of $X$ containing the diagonal and equivariant  projection-valued measures  
 $\phi^{\prime}$ and $\phi^{\prime\prime}$  on $H^{\prime}$ and $H^{\prime\prime}$
such that $(H^{\prime},\phi^{\prime})$ and $(H^{\prime\prime},\phi^{\prime\prime})$ are locally finite equivariant $X$-controlled Hilbert spaces and the inclusions of
$H^{\prime}$ and   $H^{\prime\prime}$ into $H$ have propagation controlled by   $U$.

Since $X$ is very proper we can find an   equivariant partition $(B_{i})_{i\in I}$ of $X$ with the properties listed in Definition \ref{fpowjoefewfewf} for the entourage $U$.

We assume that $I/G$ is infinite (the argument in the finite case is similar, but simpler).
Then by \ref{fpowjoefewfewf}.\ref{hweiofjweiofewf5} we can assume that $I=\bigsqcup_{n\in \nat} I_{n}$ for transitive $G$-sets $I_{n}$.
For every integer $n$ we choose a base point $i_{n}$ in $I_{n}$.

For every integer $n$ we choose a point $b_{i_{n}}$ in $B_{i_{n}}$ which is fixed by {the stabilizer} $G_{i_{n}}$ {of $i_n$}. This is possible by Assumption \ref{fpowjoefewfewf}.\ref{ojreoigjoergergerg}. For $i$ in $I_{n}$ we
then define $b_{i}:=gb_{i_{n}}$ where $g$ in $G$ is such that $gi_{n}=i$. Note that $b_{i}$ is well-defined. The collection of points $(b_{i})_{i\in I}$ thus defined is equivariant.

  We will construct the equivariant projection-valued measure $\psi$ on $H^{\prime} \mathbin{\bar{+}} H^{\prime\prime}$  by induction.

Because of the inclusion
$$\phi(B_{i_0}) H' \subseteq \phi'(V[B_{i_0}]) H^{\prime}$$   the subspace 
 $\phi(B_{i_0}) H'$ of $H$ is finite-dimensional.  Analogously we conclude that $\phi(B_{i_0}) H''$ is finite-dimensional. Consequently,   $\phi(B_{i_0}) (H' + H'')=\phi(B_{i_0}) (H' \bar{+} H'')$ is finite-dimensional.
We define the subspace $$H_{i_{0}} := \phi(B_{i_0}) (H' + H'')$$ of   $H$ and let $Q_{i_{0}}$ be the orthogonal projection onto $H_{i_{0}}$.  Note that $H_{i_{0}}$ is preserved by $G_{i_{0}}$ and that
$Q_{i_{0}}$ is $G_{i_{0}}$-invariant. We then define,  for all $i$ in $I_{0}$, the subspaces $H_{i}:=gH_{i_{0}}$ and 
$Q_{i}:= g Q_{i_{0}}g^{-1}$, where $g$ is such that $gi_{0}=i$. Note that these objects are well-defined. Furthermore, the spaces $H_{i}$ for all $i$ in $I_{0}$  are mutually orthogonal.

We let $H_{0}$ be the Hilbert subspace generated by the spaces $H_{i}$ for all $i$ in $I_{0}$.
We further define the equivariant projection-valued measure $\psi_{0}$ on $H_{0}$ by
$$\psi_{0}:= \sum_{i\in I_{0}}\delta_{b_{i}}  Q_{i}\ .$$
Using $\psi_{0}$ we recognize $H_{0}\subseteq H$  as a locally finite subspace.

Let now $n$ be an integer and assume that we have constructed an invariant subspace $H_{n}$ and an equivariant projection-valued measure $\psi_{n}$ {exhibiting the subspace $H_n$ of $H$ as a locally finite subspace.}

 As above we observe that $\phi(B_{i_{n+1}})(H'  \mathbin{\bar{+}}  H'')$ is finite-dimensional.  We define the closed subspace $$\tilde H_{i_{n+1}} := H_n + \phi(B_{i_{n+1}})(H'  \mathbin{\bar{+}}  H'') $$ of $H$ and we let  $Q_{i_{n+1}}$ be the orthogonal projection onto $H_{i_{n+1}}:=\tilde H_{i_{n+1}} \ominus H_n$.
 We note that $H_{i_{n+1}}$  and $Q_{i_{n+1}}$  are $G_{i_{n+1}}$-invariant.
For every $i$ in $I_{n+1}$ we now define
$H_{i}:=g H_{i_{n+1}}$ and $Q_{i}:=gQ_{i_{n+1}}g^{-1}$, where $g$ in $ G$ 
is such that $gi_{n+1}=i$. Note that these objects are well-defined.
Furthermore, the spaces $H_{i}$ for $i$ in $I_{n+1}$ are mutually orthogonal and orthogonal to $H_{n}$.

We let $H_{n+1}$ be the Hilbert space generated by $H_{n}$ and all the spaces $H_{i}$ for $i$ in $I_{n+1}$. This space is $G$-invariant and we can    define the equivariant projection-valued measure $\psi_{n+1}$ on $H_{n + 1}$ by $$\psi_{n+1} := \psi_n + \sum_{i\in I_{n+1}} \delta_{b_{i}}  Q_{i}\ .$$
It exhibits  $H_{n + 1} $ as a locally finite subspace of $H $.

We now observe that $$\psi:=\sum_{n\in \nat} \sum_{i\in I_{n}} \delta_{b_{i}} Q_{i}$$
is an equivariant projection-valued measure  on  the sum $H^{\prime}\bar{+}H^{\prime\prime}$, and that
the inclusion of this sum into $H$ is $V$-controlled, if $V$ is a bound for the size of the members of the partition $(B_{i})_{i\in I}$.
Furthermore, $(H^{\prime}\bar{+}H^{\prime\prime},\psi)$ is locally finite.
 \end{proof}

\section{Existence of ample \texorpdfstring{$\boldsymbol{X}$}{X}-controlled Hilbert spaces}

Let $X$ be a $G$-bornological coarse space and $(H,\phi)$ an equivariant  $X$-controlled Hilbert space.

\begin{ddd}
$(H,\phi)$ is called ample if it is determined on points and for every equivariant $X$-controlled Hilbert space $(H^{\prime},\phi^{\prime})$ there exists a controlled unitary inclusion  $H^{\prime} \to H$.
\end{ddd}

Let $X$ be a $G$-bornological coarse space.

\begin{prop}
If $X$ is very proper, then $X$ admits an ample equivariant $X$-controlled Hilbert space.
\end{prop}

\begin{proof}
Since $X$ is very proper we can find an   equivariant partition $(B_{i})_{i\in I}$ of $X$ with the properties listed in Definition \ref{fpowjoefewfewf} for some entourage $U$.

We assume that $I/G$ is infinite (the argument in the finite case is similar, but simpler).
By Assumption \ref{fpowjoefewfewf}.\ref{hweiofjweiofewf5} we can then assume that $I=\bigsqcup_{n\in \nat} I_{n}$ for a family $(I_{n})_{n\in \nat}$ of transitive $G$-sets.
For every  $n$ in $\nat$ we choose a base point $i_{n}$ in $I_{n}$.

For every $n$ in $\nat$ we choose a point $b_{i_{n}}$ in $B_{i_{n}}$ which is fixed by {the stabilizer $G_{i_{n}}$ of $i_n$.} This is possible by Assumption \ref{fpowjoefewfewf}.\ref{ojreoigjoergergerg}. For every $i$ in $I_{n}$ we
then define the point $b_{i}:=gb_{i_{n}}$ in $X$, where $g$ in $G$ is such that $gi_{n}=i$. Note that $b_{i}$ is well-defined. The collection of points $(b_{i})_{i\in I}$  is equivariant, i.e., for every $i$ in $I$ and $g$ in $G$ we have $gb_{i}=b_{gi}$.

For every {$n$} in $\nat $ we consider the $G$-Hilbert space
$$H_{n}:=(L^{2}(G)\otimes L^{2}(G_{i_{n}}))^{G_{i_{n}}}\otimes \ell^{2}\ ,$$
where  the $G_{i_{n}}$-invariants are taken with respect to the action of $G_{i_{n}}$ on $L^{2}(G)$ by right-multiplication on $G$ and on $L^{2}(G_{i_{n}})$  by the left-multiplication on $G_{i_{n}}$. The $G$-action on $H_{n}$ is induced from the left-multiplication of $G$ on itself.
We then define the $G$-Hilbert space
$$H:=\bigoplus_{n\in \nat} H_{n}\ .$$ 
We interpret the elements of $H_{n}$ as functions on $G$ with values in $L^{2}(G_{i_{n}})\otimes \ell^{2}$.
For every $n$ in $\nat$  and $i$ in $I_{n}$ we let
$Q_{i}$ be the projection onto
the subspace $H_{i}$ of functions in $H_{n}$ supported on the subset $\{g\in G\:|\: gi_{n}=i\}$ of $G$.   Since this subset is preserved by the right-action $G_{i_{n}}$ the support condition is compatible with the $G_{i_{n}}$-invariance condition, and $H$ is spanned by these subspaces  $H_{i}$ for all  $i$ in $I$.    We further note that $gQ_{i}g^{-1}= Q_{gi}$.
  
  We define the measure
$$\phi:=\sum_{i\in I} \delta_{b_{i}}Q_{i}$$
which turns out to be invariant.
In this way we get an equivariant $X$-controlled Hilbert space $(H,\phi)$.
It is determined on points.

We now show that $(H,\phi)$ is ample.
Let $(H^{\prime},\phi^{\prime})$ be any equivariant $X$-controlled Hilbert space.
Then we must construct a controlled  isometric embedding $H^{\prime}\to H$.

For every $i$ in $I$ we consider the subspace $H_{i}^{\prime}:=H^{\prime}(B_{i})$. Since the subset $B_{i}$ is bounded and $G_{i}$-invariant, 
this is a separable representation of the finite {(by Assumption \ref{fpowjoefewfewf}.\ref{hweiofjweiofewf})} group $G_{i}$.

The action of $G$ on $H$ restricts to an action   of $G_{i}$ on $H_{i}$
{which is} of the form $L^{2}(G_{i})\otimes \ell^{2}$. By the Peter--Weyl theorem it contains
every {irreducible} representation of $G_{i}$ with countably infinite multiplicity. 
{Hence} for every $n$ we can choose an isometric embedding
$u_{i_{n}}:H_{i_{n}}^{\prime}\to H_{i_{n}}$ which is $G_{i_{n}}$-invariant.
We extend this $G$-equivariantly by
$$u_{i}:=gu_{i_{n}}g^{-1}:H_{i}^{\prime}\to H_{i}\ ,$$
where $g$ is any element of $G$ such that $gi_{n}=i$.
The sum $u:=\oplus_{i\in I} u_{i}$ is then an equivariant embedding of $H^{\prime}$ into $H$.
By Assumption \ref{fpowjoefewfewf}.\ref{iegergerg} there exists an  entourage $V$ of $X$  which bounds the size of the members of the partition.
The embedding  $u$ is then $V$-controlled.
\end{proof}

\begin{rem}
The above proof would also work without Assumption \ref{fpowjoefewfewf}.\ref{hweiofjweiofewf5} that $I/G$ is countable. Moreover we have  not used Assumption \ref{fpowjoefewfewf}.\ref{hweiofjweiofewf1}. Furthermore, for this construction we can also replace finite by countable  
 in Assumption \ref{fpowjoefewfewf}.\ref{gergihieorgergrege}. This still ensures that $H(B)$ is separable for every bounded subset $B$ of $X$.
\hB
\end{rem}

\section{\texorpdfstring{$\boldsymbol{C^{*}}$}{C*}-categories and the construction of \texorpdfstring{$\boldsymbol{\KX^{G}}$}{KXG}}\label{kfopkwefewdqwedcacsdcac}

In this section we introduce the $C^{*}$-categories of locally finite equivariant controlled Hilbert spaces and define the equivariant coarse $K$-homology functor.
This material is taken from \cite{bu} and  reproduced here for the sake of self-containedness.

\begin{ddd}
We let $\bC(X)$ 
be the $\C$-linear $*$-category given by the following data: 
\begin{enumerate}
\item The objects of $\bC(X)$ are equivariant locally finite  $X$-controlled Hilbert spaces.
\item The $\C$-vector space of morphisms between the objects $(H^{\prime},\phi^{\prime})$ and  $(H,\phi)$ is defined as $C_{\lc}(X,(H^{\prime},\phi^{\prime}),(H,\phi))$. 
\item The involution $*$ is given by taking adjoints.
\end{enumerate}
\end{ddd}

Note that any operator between locally finite  $X$-controlled Hilbert spaces is locally compact. As explained in \cite{Bunke:2016aa}, being a $C^{*}$-category is a property of a $\C$-linear $*$-category.
 One can check that
$\bC(X)$ 
is a  $C^{*}$-category.

If $f:X\to X^{\prime}$ is a morphism of $G$-bornological coarse spaces, then we get a functor between $C^{*}$-categories
$$f_{*}:\bC(X)\to \bC(X^{\prime})
$$ given as follows:
\begin{enumerate}
\item $f_{*}$ sends the object $(H,\phi)$ to $(H,f_{*}\phi)$, where $(f_{*}\phi)(Y'):=\phi(f^{-1}(Y))$ for every subset $Y'$ of $X'$.
\item $f_{*}$ is given by the identity on morphisms.
\end{enumerate}
  This defines   a functor
\begin{equation}\label{huihwiefewrfwfwrfwerfw}\bC
:G\BC\to\Ccat\ ,
\end{equation} 
where $\Ccat$ denotes the category of unital $C^{*}$-categories and functors.
We let
 $$K:\Ccat\to \Sp$$ be a topological $K$-theory functor from $C^{*}$-categories to spectra (a discussion of constructions of such a functor may be found in \cite[Sec.~8.4 \& 8.5]{buen}). 
 The $K$-theory functor for  unital $C^{*}$-categories used above is the restriction of a more general 
functor for non-unital $C^{*}$-categories
$K:\nCcat\to \Sp$, see \cite[Def. 14.4]{cank}.  
 Note that a $C^{*}$-algebra $A$ can be considered as a possibly non-unital $C^{*}$-category with one object, and $K(A)$ is the usual usual $K$-theory spectrum of $A$ which is essentially uniquely characterized
 by \cite[Prop. 2.7]{blp}
 

\begin{ddd}\label{girogegegergerg}
We define the equivariant coarse $K$-homology functor
$$\KX^{G}:G\BC\to \Sp$$
by 
$$\KX^{G}:=K\circ \bC\ .$$
\end{ddd}

\begin{theorem}[\cite{bu}]
$\KX^{G}$ 
is an equivariant coarse homology theory.
\end{theorem}

In this paper we will not recall the notion of an equivariant coarse homology theory but rather refer to \cite{equicoarse} for definitions.
But we remark that $\KX^{G}$ in particular sends coarse equivalences to equivalences of spectra and satisfies excision for invariant, coarsely excisive decompositions.

\section{Comparison with the \texorpdfstring{$\boldsymbol{K}$}{K}-theory of Roe algebras}

Let  $X$ be a $G$-bornological coarse space and let
$(H,\phi)$ be an equivariant $X$-controlled Hilbert space.
Let $X^{\prime}$ and  $(H^{\prime},\phi^{\prime})$ be similar data. \begin{theorem}\label{fwefiwjfeiooi234234324434e}
\mbox{}
\begin{enumerate}
\item Assume that $X$ is very proper  {and that $(H,\phi)$ is ample.}
 Then there exists a canonical (up to equivalence) equivalence of spectra  
$$\kappa_{(X,H,\phi)}:K(C(X,H,\phi))\to \KX^{G}(X)\ .$$
\item Assume that $X$ and $X^{\prime}$ are very proper and that $(H,\phi)$ and $(H^\prime,\phi^\prime)$ are ample. For a morphism  $f:X^{\prime}\to X$  of $G$-bornological coarse spaces and an equivariant isometric inclusion $V:H^{\prime}\to H$   which induces a controlled morphism $f_{*}(H^{\prime},\phi^{\prime})\to (H,\phi)$ we have a commutative diagram
\begin{equation}\label{coiehjocwecwecc}
 \xymatrix{K(C(X^{\prime},H^{\prime},\phi^{\prime})) \ar[rr]^(0.55){\kappa_{(X^{\prime},H^{\prime},\phi^{\prime})}}\ar[d]^{ K(v)}&& \KX^{G}(X^{\prime})\ar[d]^{f_{*}}\\
 K(C(X,H,\phi))   \ar[rr]^(0.55){\kappa_{(X,H,\phi)}}&& \KX^{G}(X)}
\end{equation}
\end{enumerate}
\end{theorem}

\begin{proof}
Given Proposition \ref{iofeoewfewfewf} the proof is analogous to the proof of \cite[Thm.~8.99]{buen}.
Since we need some of the details later we will recall the construction of the equivalence $\kappa_{(X,H,\phi)}$.

We consider the category $\bC(X,H,\phi)$ whose objects are triples 
$(H^{\prime},\phi^{\prime},U)$, where   $(H^{\prime},\phi^{\prime})$ is an object of $\bC(X)$, and $U$ is an equivariant  controlled inclusion
$U:H^{\prime}\to H$ as a locally finite subspace. We consider  the Roe algebra $C(X,H,\phi)$ as a non-unital  $C^{*}$-category   with a single object.
Then we have  functors between non-unital $C^{*}$-categories
\begin{equation}\label{rgejknerkggergegergrege}
\bC(X) \stackrel{F}{\leftarrow}  \bC(X,H,\phi) \stackrel{I}{\to} C(X,H,\phi)\ .
\end{equation} 
The functor $F$ forgets the inclusions. The action of $I$ on objects is clear. Furthermore it maps the morphism $A:(H^{\prime},\phi^{\prime},U^{\prime})\to  (H^{\prime\prime},\phi^{\prime\prime},U^{\prime\prime})$  in   $\bC(X,H,\phi)$ to the morphism   $I(A):=U^{\prime\prime}AU^{\prime,*}$  of  $C(X,H,\phi)$. 

We get an induced diagram
of $K$-theory spectra 
\begin{equation}\label{rgejknerkggergegergrege1}
K( \bC(X))  \stackrel{K( F) }{\leftarrow} K(    \bC(X,H,\phi) ) \stackrel{K( I)  }{\to} K(C(X,H,\phi))\ .\end{equation} 

Using Proposition \ref{iofeoewfewfewf}, as in the  proof    of \cite[Thm.~8.99]{buen} we now see that $K(F)$ and $K(I)$ are equivalences.  We set
$\kappa_{(X,H,\phi)}:=K(F)\circ K(I)^{-1}$.

The proof of the second assertion is the same as in  \cite[Thm.~8.99]{buen}.
\end{proof}

Combining Theorem \ref{fwefiwjfeiooi234234324434e} with Proposition \ref{lem:sdfbi23} we immediately get the following corollary.

Let $X$ be a $G$-bornological coarse space and let $(H,\phi)$ be an equivariant $X$-controlled Hilbert space.
\begin{kor}
If $X$ is very proper and $(H,\phi)$ is ample, then we have a canonical (up to equivalence) equivalence of spectra $K(C_{\lc}(X,H,\phi)) \simeq \KX^G(X)$.
\end{kor}

\section{\texorpdfstring{$\boldsymbol{K}$}{K}-theory of graded \texorpdfstring{$\boldsymbol{C^{*}}$}{Cstar}-algebras and \texorpdfstring{$\boldsymbol{C^{*}}$}{Cstar}-categories}\label{erkogpqwerfrewfrefrefw}

\newcommand{\RoeCat}{\overline{\bD}^{G,\mathrm{ctr}}_{\mathrm{lf}}}
\newcommand{ \grRoeCat}{{}^{\mathrm{gr}}\overline{\bD}^{G,\mathrm{ctr}}_{\mathrm{lf}}}
\newcommand{\EE}{\mathrm{E}}
\newcommand{\ee}{\mathrm{e}}
 \renewcommand{\bG}{\mathbf{G}}
  \newcommand{\bH}{\mathbf{H}}
  \newcommand{\triv}{\mathrm{triv}}

In this section we recall some basic facts about the $K$-theory of graded $C^{*}$-algebras and $C^{*}$-categories. We furthermore extend the construction of coarse equivariant $K$-homology to some graded coefficient categories and extend Theorem \ref{fwefiwjfeiooi234234324434e}. The main motivation is to capture the index of a Dirac-type operator  which is graded and commutes with the action of some Clifford algebra.

 The category of graded $C^{*}$-algebras $\grnCalg$
is equivalent to the category $C_{2}\nCalg$ of $C^{*}$-algebras with $C_{2}$-action and equivariant maps, where 
 $C_{2}$ denotes the group with two elements.   But in contrast to the latter the tensor products $\hat \otimes_{\min}$ and $\hat \otimes_{\max}$ in $\grnCalg$ will be defined using Koszul sign rules. Furthermore, the $K$-theory of graded algebras differs from the $C_{2}$-equivariant $K$-theory.    Following \cite{MR1694805}, \cite{Meyer:aa}, \cite{zbMATH02068495}  a quick way to define it as a spectrum-valued functor is, using equivariant $E$-theory, as 
\begin{equation}\label{ewfwedqwededq}\hat K \colon \grnCalg\to \Sp\ , \quad \hat K(-) \coloneqq \EE^{C_{2}}(\hat C_{0}(\R),-)\ ,  
\end{equation}  where  $\hat C_{0}(\R)$ has the $C_{2}$-action induced by multiplication by $-1$ on $\R$.
Here we use the equivariant $E$-theory functor $\ee^{C_{2}} \colon C_{2}\nCalg\to \EE^{C_{2}}$ from \cite[Sec.~3.2]{budu}
(the universal homotopy invariant, $K_{C_{2}}$-stable, exact, s-finitary and countable sum-preserving functor to a cocomplete stable $\infty$-category) and the notation
$$\EE^{C_{2}}(A,B):=\map_{\EE^{C_{2}}}(\ee^{C_{2}}(A),\ee^{C_{2}}(B))$$
for the morphism spectrum. 
  In contrast to the references   listed above we prefer $E$-theory over $K\!K$-theory since this makes the exactness property of the $K$-theory functor transparent. 
\begin{ex}
 Let $\Cl^{n}$ denote the unital graded complex Clifford algebra of the Euclidean space $\R^{n}$. It is generated by the linear subspace $\R^{n}$ consisting of 
 anti-selfadjoint elements $x$  subject to the relation $x^{2}= -\|x\|^{2}$. The grading is induced by
 functoriality and the multiplication by $-1$ on $\R^{n}$. In particular, the generators $x$ are odd.

 For an arbitrary graded $C^{*}$-algebra $A$ 
 we have natural equivalences $ \hat K(A)\simeq  \Sigma^{n} \hat K(A\hat \otimes \Cl^{n})$ inducing isomorphisms
$ \hat K_{*}(A)\cong  \hat K_{*-n}(A\hat \otimes \Cl^{n})$ on the level of $K$-groups, see Example \ref{erigjorewgwrferfrwefwref}.  
\hB
 \end{ex}

\begin{prop}\label{jiefgoerfqewefqewd}The $K$-theory for graded algebras extends the one for ungraded algebras: We have a commutative diagram
$$\xymatrix{\nCalg\ar[rr]^{K}\ar[dr]_{\incl} && \Sp  \\ &\grnCalg\ar[ur]_{\hat K} & } $$
of functors from $\nCalg$ to $\Sp$.
\end{prop}
\begin{proof}
 Let $\Res_{C_{2}}:\nCalg\to C_{2}\nCalg$ be the functor which equips a
 $C^{*}$-algebra with the trivial $C_{2}$-action. If we interpret the target as graded $C^{*}$-algebras, then this functor is the inclusion featuring in the statement.  In view of  \eqref{ewfwedqwededq} and the equivalence $K(-)\simeq \EE(\C,-)$ we must construct an equivalence of functors 
 $$\EE^{C_{2}}( \hat C_{0}(\R),\Res_{C_{2}}(-))\simeq  \EE(\C,-):\nCalg\to \Sp\ .$$
The evaluation at $0$ induces the   morphism $\epsilon$ in an exact sequence
\begin{equation}\label{frewfwerfewrfws}0\to S(\Cl^{1})\to \hat C_{0}(\R)\xrightarrow{\epsilon} \C\to 0
\end{equation} 
in $C_{2}\nCalg$. In order to describe the first morphism we identify $\Cl^{1}\cong C(\{-1,1\})$, where   the $C_{2}$-action on the right-hand side exchanges $1$ with $-1$.   The  first map then sends
$f$ in $S(\Cl^{1})\cong C_{0}((0,\infty) ,\Cl^{1})$ to the function
$t\mapsto f(t^{2})(\sign(t))$. 

Let $\ee : \nCalg\to \EE$ denote the non-equivariant $E$-theory functor.
 The two one-dimensional representations $\triv$ and $\sign$ of $C_{2}$ provide  
homomorphisms $\C\rtimes C_{2}\to \C$ which induce an equivalence 
\begin{equation}\label{qewfqweddedqed}(\triv,\sign):\ee(\C\rtimes C_{2})\simeq \ee(\C)\oplus \ee(\C)
\end{equation} 
in the stable $\infty$-category $\EE$.
The composition 
\begin{equation}\label{frewfrefwrefwfref}\ee( \hat C_{0}(\R)\rtimes C_{2})\xrightarrow{\ee(\epsilon\rtimes C_{2})} \ee(\C\rtimes C_{2})\xrightarrow{\sign}
 \ee(\C)
\end{equation} is an equivalence.

 By the dual Green--Julg theorem  \cite[Prop. 3.61.2]{budu}
 we have an equivalence
 $$\EE^{C_{2}}( \hat C_{0}(\R),\Res_{C_{2}}(-))\simeq \EE( \hat C_{0}(\R)\rtimes C_{2},-):\nCalg\to \Sp\ .$$
 Inserting the equivalence \eqref{frewfrefwrefwfref} in the first entry we get the desired equivalence
\[\EE^{C_{2}}( \hat C_{0}(\R),\Res_{C_{2}}(-))\simeq  \EE(\C,-) \ .\qedhere\]
\end{proof}

\begin{prop}\label{erkogpwergfefrefwf}
We have a fibre sequence  \begin{equation}\label{fqewfweewdqdewd}
   K(-\rtimes C_{2})\to     \hat K(-)\to \Omega K(\Res^{C_{2}}(-))
\end{equation}
of functors from $\grnCalg$ to $\Sp$.
\end{prop}
\begin{proof}
We use the 
equivalence
$\Cl^{1}\simeq \Ind^{C_{2}}(\C)$ as a $C_{2}$-$C^{*}$-algebra,
the $(\Ind^{C_{2}},\Res^{C_{2}})$-adjunction on the level of $E$-theory (see \cite[Prop. 3.60]{budu}), and the Green--Julg  equivalence $\EE^{C_{2}}(\C,-)\simeq 
\EE(\C,-\rtimes C_{2})$. 
 Plugging the exact sequence \eqref{frewfwerfewrfws} into the first entry of $\EE^{C_{2}}$
and using exactness of $\ee$ and the equivalences $\EE^{C_{2}}(\C,-)\simeq K(-\rtimes C_{2})$ 
and $\EE^{C_{2}}(S(\Cl^{1}),-)\simeq \Omega K(\Res^{C_{2}}(-))$  we get the desired  fibre sequence \end{proof}

%
%

   From now one we will omit the superscript and just write $K(A)$ also for the $K$-theory of a graded algebra.

 We next extend the $K$-theory of graded $C^{*}$-algebras to graded $C^{*}$-categories.
A graded $C^{*}$-category is a $C^{*}$-category with a strict $C_{2}$-action
fixing objects. A morphism of graded $C^{*}$-categories is an equivariant functor.
 Therefore the graded $C^{*}$-categories form a full subcategory $\grnCcat$ of the category $C_{2}\nCcat$ of $C^{*}$-categories with a strict $C_{2}$-action.

\begin{ex}
Our typical example of a graded $C^{*}$-category  is the category  $\Hilb_{c,(n)}(\C)$ of graded Hilbert spaces  with a right action of $\Cl^{n}$ and compact operators which commute with $\Cl^{n}$.

A graded $C^{*}$-algebra can be viewed as a graded $C^{*}$-category with a single object. We thus get an inclusion functor  $\incl:\grnCalg\to \grnCcat$.
\hB  \end{ex}

By functoriality the adjunction  $$A^{f}:\nCcat \leftrightarrows \nCcat:\incl$$ (see \cite[Def. 3.7]{joachimcat}) induces a functor
$$A^{f}:\grnCcat\leftrightarrows \grnCalg:\incl\ .$$ Generalizing \cite[Def. 14.3]{cank} from the ungraded to the graded case   we define the $K$-theory functor for graded $C^{*}$-categories as the composition
$$K:\grnCcat\xrightarrow{A^{f}}  \grnCalg\xrightarrow{K} \Sp\ .$$

 Proposition \ref{jiefgoerfqewefqewd} has the following consequence.
 \begin{kor}The $K$-theory for graded  $C^{*}$-categories extends the one for ungraded $C^{*}$-categories. \end{kor}
 \begin{lem}\label{wergkoepferfwerf}
 We have a fibre sequence  \begin{equation}\label{fqewfweewdqdewd76}
   K(-\rtimes C_{2})\to      K(-)\to \Omega K(\Res^{C_{2}}(-))\end{equation} of functor from $\grnCcat$ to $\Sp$.
\end{lem}
\begin{proof}
This follows from Proposition \ref{erkogpwergfefrefwf}, the equivalence
$ \ee(A^{f}(-\rtimes C_{2}))\simeq \ee(A^{f}(-)\rtimes C_{2})$, and the isomorphism  $A^{f}(\Res_{C_{2}}(-))\simeq \Res^{C_{2}}(A^{f}(-))$.
\end{proof}

 An exact sequence $$0\to \bC\to \bD\to \bE\to 0$$ of graded $C^{*}$-categories is a sequence $\bC\to \bD\to \bE$ of morphisms between graded $C^{*}$-categories which becomes an exact sequence after forgetting the grading.
Recall that we can form the crossed product of $C_{2}$ with a $C^{*}$-category with $C_{2}$-action, and that this  crossed product preserves exact sequences \cite{crosscat}.

A morphism $\bC\to \bD$ in  $\grnCcat$  is a relative Morita equivalence if there exists a diagram
$$
\xymatrix{
0\ar[r]&\bC\ar[d]\ar[r]&\bE\ar[r]\ar[d]&\bF\ar[r]\ar[d]&0\\
0\ar[r]&\bD\ar[r]&\bG\ar[r]&\bH\ar[r]&0}
$$
of exact sequence in $\grnCcat$ where $\bE,\bF,\bG,\bH$ are unital and
$\bE\to \bG$ and $\bF\to \bH$ are Morita equivalences  of unital $C^{*}$-categories after forgetting 
the grading \cite[Def. 17.1]{cank}.

 \begin{prop}\label{wkogrperfrefw}
  The functor $K:\grnCcat\to \Sp$  has the following properties: 
  
  \begin{enumerate}
\item  $K$ sends exact sequences to fibre sequences.
\item  $K $ preserves filtered colimits.
\item $K$ sends relative Morita equivalences to equivalences.
\end{enumerate}
 \end{prop} \begin{proof}By Lemma \ref{wergkoepferfwerf}
 we have an equivalence $$K(-)\simeq \Fib( \Omega K(\Res^{C_{2}}(-))\to   K(-\rtimes C_{2}))\ .$$
 Since the two functors $ \Omega K(\Res^{C_{2}}(-))$ and $   K(-\rtimes C_{2})$ 
 send exact sequences to fibre sequences, preserve filtered colimits, and send  relative Morita equivalences to equivalences so does $K$.
%
%
%
%
\end{proof}

%

We say that a graded $C^{*}$-category $\bC$   admits countable AV-sums or  is effectively additive if the underlying $C^{*}$-category does so (see \cite[Sec. 3]{bu} for these notions).  

\begin{ex}
The category  $\Hilb_{c,(n)}(\C)$    is effectively additive and admits countable AV-sums.   \hB  \end{ex}

We next generalize the construction of the equivariant coarse $K$-homology $K\cX^{G}_{\bD}$ with coefficients in a $C^{*}$-category   $\bD$   with strict $G$-action \cite[Thm. 7.3]{bu} to the graded case. Thus
 let $\bD$ be in $\Fun(BG,\grnCcat)$ be a graded $C^{*}$-category with a strict $G$-action.
 As in the ungraded case we assume that   it admits countable AV-sums and  is effectively additive. Recall that the objects
 of  the $C^{*}$-category $\RoeCat(X)$ of $G$-equivariant $X$-controlled objects  in $\bD$   defined in \cite[Def. 5.24]{bu} are  triples $(D,\rho,\mu)$ of an object $D$ of $\bD$, a cocycle of unitaries  $\rho=(\rho_{g}:D\to gD)$, and a projection-valued finitely-additive measure $\mu$. We let  $\grRoeCat(X)$ denote the full subcategory of  $\RoeCat$ on objects $(D,\rho,\mu)$ where $\rho$ and $\mu$ consist of even operators. In this way we get a subfunctor $\grRoeCat$ of $\RoeCat$.

  \begin{theorem} The functor
  $$K\cX_{\bD}^{G}:G\BC\xrightarrow{\grRoeCat} \grnCcat \xrightarrow{K} \Sp\ .$$
  is a coarse homology theory.
\end{theorem}
\begin{proof}
The proof is the same as \cite[Thm. 7.3]{bu}.
It also uses Proposition \ref{wkogrperfrefw} saying that the $K$-theory functor for graded $C^{*}$-categories has the same exactness and continuity  properties as the usual one for ungraded categories.
\end{proof}

 Let  $n$ be in $\nat$. 
 In the following we make the construction of $\bC_{(n)}:=\grRoeCat$ for $\bD=\Hilb_{c,(n)}(\C)$ explicit.
 Let $X$ be a $G$-bornological coarse space.

\begin{ddd}An equivariant $X$-controlled Hilbert space of degree $n$ is an equivariant  $X$-controlled Hilbert space $(H,\phi)$ with an additional grading and a graded action of $\Cl^{n}$.  We require that the measure $\phi$ commutes with the action of $\Cl^{n}$ and preserves the grading.
\end{ddd}

 We use the notation $(H,\phi,n)$ in order to denote an  equivariant $X$-controlled Hilbert space of degree $n$.

Let $(H,\phi,n)$ be an equivariant  $X$-controlled Hilbert space of degree $n$.
\begin{ddd} We let 
$C_{\lc}(X,H,\phi,n)$ be the graded subalgebra of $C_{\lc}(X,H,\phi)$   of operators in $C_{\lc}(X,H,\phi)$ which commute with the action of $\Cl^{n}$.
\end{ddd}

The definitions above translate to the following:
\begin{prop}
 The  graded  $C^{*}$-category $\bC_{(n)}(X)$ 
is  given by the following data: 
\begin{enumerate}
\item The objects of $\bC_{(n)}(X)$ are equivariant locally finite  $X$-controlled Hilbert spaces $(H,\phi,n)$ of degree $n$.
\item The  $\C$-vector space of morphisms between the objects $(H^{\prime},\phi^{\prime},n)$ and  $(H,\phi,n)$ is  $C_{\lc}(X,(H^{\prime},\phi^{\prime},n),(H,\phi,n))$ with its canonical grading. 
\item The involution is given by taking adjoints.
\end{enumerate}
\end{prop}

For $n$ in $\nat$ we use the notation  $$K\cX^{G}_{(n)}:=K\cX^{G}_{\Hilb_{c,(n)}(\C)}=K\circ \bC_{(n)}$$
in order to denote the equivariant coarse $K$-homology theory associated to the coefficient category graded $C^{*}$-category $ \Hilb_{c,(n)}(\C)$.

\begin{prop}\label{kgowpegrfwfwfref}
We have an equivalence of equivariant coarse homology theories
$$K\cX^{G}_{(n)}\simeq \Omega^{n} K\cX^{G}\ .$$
\end{prop}
\begin{proof}
We use obvious extension of the tensor product of $C^{*}$-categories from the ungraded to the graded case (see \cite[Sec. 7]{KKG} for a detailed discussion).  
We consider $\Cl^{1}$ as a graded $C^{*}$-category with a single object so that $-\hat \otimes \Cl^{1}$ has the effect of forming the graded tensor product of all morphism spaces with $\Cl^{1}$.
We have a natural fully faithful functor 
\begin{equation}\label{dtd78qt7zd8qwswqs} \bC_{(n)}(X)\hat \otimes \Cl^{1}\stackrel{\simeq}{\to} \bC_{(n+1)}(X)
\end{equation}    sending the object $(H,\phi,n)$ to $ (H\hat \otimes \Cl^{1},\phi\hat \otimes \id_{\Cl^{1}},n+1)$, and a morphism
  $A\hat \otimes \sigma$ in $$\Hom_{ \bC_{(n)}(X)}((H,\phi,n),(H',\phi',n))\hat\otimes \Cl^{1} $$ to itself considered as a morphism $H\hat\otimes \Cl^{1}\to H'\hat\otimes \Cl^{1}$, where $\sigma$ acts by left multiplication.

We claim that the functor in \eqref{dtd78qt7zd8qwswqs} is a relative Morita equivalence.
Indeed, if $n+1$ is even, then  $\Cl^{n+1} $
 as ungraded algebra has a unique irreducible representation $S_{n+1}$ of dimension $2^{\frac{n+1}{2}}$. 
 Hence after forgetting the grading every object of $ \bC_{(n+1)}(X)$ is isomorphic to an object of the form
 $(H'\hat \otimes S_{n+1},\phi'\hat \otimes\id_{S_{n+1}},n+1)$, and the same is true for
 $\bC_{(n)}(X)\hat  \otimes \Cl^{1}$.
 If $n+1$ is odd, then $\bC_{(n+1)}(X)$ has two ungraded representations $S^{\pm}_{n}$ of dimension  $2^{\frac{n}{2}}$. Every object of $\bC_{(n+1)}(X)$
 is isomorphic to a sum   $$(H'_{+}\hat \otimes S^{+}_{n+1},\phi'_{+}\hat \otimes\id_{S^{+}_{n+1}},n+1)\oplus (H'_{-}\hat \otimes S^{-}_{n+1},\phi'_{-}\hat \otimes\id_{S^{-}_{n+1}},n+1)\ .$$
After forgetting the grading the objects in the tensor product  $\bC_{(n)}(X)\hat \otimes \Cl^{1}$ are isomorphic to such sums with $H'_{+}\cong H'_{-}$ and $\phi_{+}=\phi_{-}$ (under this isomorphism).
We can conclude that every object of $\bC_{(n)}(X)$ is a summand of an object in the image of the functor.
 
  We therefore have equivalences of equivariant coarse homology theories 
$$K\cX^{G}_{(n)}\simeq \Sigma^{}K\cX^{G}_{(n-1)}\simeq \dots\simeq \Sigma^{n}K\cX^{G}_{(0)}\ .$$
It remains to show that
the inclusion $\bC(X)\to \bC_{(0)}(X)$ induces an equivalence in $K$-theory for every $G$-bornological coarse space. But this inclusion
is a Morita equivalence. In fact, after forgetting the grading it is a unitary eqivalence. We conclude that
$K\cX^{G}\simeq K\cX_{(0)}^{G}$.
  \end{proof}



\begin{ddd}
An equivariant $X$-controlled Hilbert space  $(H,\phi,n) $   of degree $n$ is ample if it is determined on points and for every 
equivariant $X$-controlled Hilbert space $(H^{\prime},\phi^{\prime},n)$  of degree $n$
there exists a controlled isometric embedding $H^{\prime}\to H$ which preserves the grading and  is $\Cl^{n}$-linear.
\end{ddd}

If $(H,\phi,n)$ is an ample equivariant $X$-controlled Hilbert space of degree $n$ and $(H^{\prime},\phi^{\prime},n)$ is an arbitrary equivariant $X$-controlled Hilbert space of degree $n$, then
we can  find a $\Cl^{n}$-linear controlled isometric embedding
$H^{\prime}\oplus H\to H$. In this case we say that the induced embedding $H^{\prime}\to H$ has ample complement.

We consider an equivariant $X$-controlled Hilbert space  $(H,\phi,n) $ of degree $n$ and  any other 
equivariant $X$-controlled Hilbert space $(H^{\prime},\phi^{\prime},n)$  of degree $n$.

\begin{lem}\label{irjfiorregregregre}
Any two controlled isometric embeddings $(H^{\prime},\phi^{\prime},n)\to (H,\phi,n)$  with the property that the complement of the sum of their images is ample are homotopic to each other.
\end{lem}

\begin{proof}
Let $U$ and $V$ by two such embeddings.
Assume {first} that the images are orthogonal. Then
$\cos(t) U+\sin(t) V$  for $t$ in $[0,\pi/2]$ is a homotopy between $U$ and $V$.

We now consider the general case. {We can choose a controlled isometric embedding $W$ of $H\oplus H$ into the complement of the sum of the images of the embeddings $U$ and $V$.} Composing $W$ with the embedding $H\to H\oplus H$, $x\mapsto x\oplus 0$ and $V$ we get an embedding $V^{\prime}$.  Similarly, composing $W$ with the embedding $H\to H\oplus H$, $x\mapsto 0\oplus x$ and $U$ we get an embedding $U^{\prime}$. 
By the above $U$ and $U^{\prime}$, $U^{\prime}$ and $V^{\prime}$, and $V^{\prime}$ and $V$ are homotopic to each other.
\end{proof}
 
Let $(H,\phi)$ be an ample equivariant $X$-controlled Hilbert space. Then we form the graded Hilbert space  $$\hat H:=(H\oplus H^{op})\hat\otimes \Cl^{n}$$ with the projection valued measure $$\hat \phi:=(\phi\oplus \phi)\hat\otimes \id_{\Cl^{n}}\ .$$
It has an $\Cl^{n}$-action   induced from the  right-multiplication on the $\Cl^{n}$-factor.  We get an equivariant $X$-controlled Hilbert space $(\hat H,\hat \phi ,n)$ of degree $n$.
 
We note that the construction of $\hat H$ implicitly depends on the degree though this is not reflected by the notation. 

\begin{lem} \label{okgpergergreg} $(\hat H,\hat \phi ,n)$   is ample. \end{lem} 

\begin{proof}
Let $(H^{\prime},\phi^{\prime},n)$ be a second  equivariant $X$-controlled Hilbert space of degree $n$.
Then we can choose an even equivariant controlled  {isometric} embedding
$i:H^{\prime}\to H\oplus H^{op}$. We then define, using a standard  basis $(e_{i})_{i=1}^{n}$ of $\R^{n}$, 
$$U:H^{\prime}\to  {(H\oplus H^{op})}\hat \otimes \Cl^{n}\ , \quad U(h):=\sum_{(k,i_{1}<\dots<i_{k})} (-1)^{k}i(he_{i_{k}}\cdots e_{i_{1}})\hat\otimes e_{i_{1}} \cdots e_{i_{k}}\ .$$
This map is a controlled,  equivariant and $\Cl^{n}$-equivariant {isometric} embedding (for a suitably normalized scalar product on $\Cl^{n}$).
\end{proof}

We now generalize Theorem \ref{fwefiwjfeiooi234234324434e} to the graded case.
Let  $X$ be a $G$-bornological coarse space and let
$(H,\phi)$ be an equivariant $X$-controlled Hilbert space. 
Let $X^{\prime}$ and  $(H^{\prime},\phi^{\prime})$ be similar data.
We use the notation explained before Lemma \ref{okgpergergreg}.
 \begin{theorem}\label{fwefiwjfeiooi234234324434en} \mbox{}
\begin{enumerate}
\item Assume that $X$ is very proper  {and that $(H,\phi)$ is ample.}
 Then there exists a canonical (up to equivalence) equivalence of spectra  
$$\kappa_{(X,H,\phi,n)}:K(C(X,\hat H,\hat \phi,n))\to \KX_{(n)}^{G}(X)\ .$$
\item Assume that $X$ and $X^{\prime}$ are very proper and that $(H,\phi)$ and $(H^\prime,\phi^\prime)$ are ample. For a morphism  $f:X^{\prime}\to X$  of $G$-bornological coarse spaces and an equivariant isometric inclusion $V:H^{\prime}\to H$   which induces a controlled morphism $f_{*}(\hat H^{\prime},\hat \phi^{\prime},n)\to (\hat H,\hat \phi,n)$ we have a commutative diagram
\begin{equation}\label{coiehjocwecwecc78}
 \xymatrix{K(C(X^{\prime},\hat H^{\prime},\hat \phi^{\prime},n)) \ar[rr]^(0.55){\kappa_{(X^{\prime},\hat H^{\prime},\hat \phi^{\prime},n)}}\ar[d]^{ K(v)}&& \KX_{(n)}^{G}(X^{\prime})\ar[d]^{f_{*}}\\
 K(C(X,\hat H,\hat \phi,n))   \ar[rr]^(0.55){\kappa_{(X,\hat H,\hat \phi,n)}}&& \KX_{(n)}^{G}(X)}
\end{equation}
\end{enumerate}
\end{theorem}
\begin{proof}
The proof is  essentially the same as for   Theorem \ref{fwefiwjfeiooi234234324434e} and based on a zig-zag
\begin{equation}\label{rgejknerkggergegergre2ge}
\bC_{(n)}(X) \stackrel{F_{(n)}}{\leftarrow}  \bC_{(n)}(X,\hat H,\hat \phi) \stackrel{I_{(n)}}{\to} C(X,\hat H,\hat \phi,n)\ .
\end{equation} 
Here $ \bC_{(n)}(X,\hat H,\hat \phi)$ consists of pairs $((H',\phi',n),U)$ of objects of $\bC_{(n)}(X)$ and an even controlled $\Cl^{n}$-linear  isometric embedding $U:(H',\phi',n)\to (\hat H,\hat \phi,n)$.
The functor $F_{(n)}$ forgets $U$ and is a unitary equivalence.  It therefore induces an equivalence in $K$-theory. We must show that $I_{(n)}$ induces an equivalence in $K$-theory. To this end first consider the square
$$\xymatrix{  \bC_{(n+1)}(X,\hat H,\hat \phi) \ar[r]^{I_{(n+1)}}& C(X,\hat H,\hat \phi,n+1)  \\ \ar[u]\bC_{(n)}(X,\hat H,\hat \phi) \hat\otimes \Cl^{1} \ar[r]^{I_{n}\hat\otimes \Cl^{1}} &C(X,\hat H,\hat \phi,n)\hat\otimes \Cl^{1}\ar[u]^{\cong} }$$
 Since the left vertical map induces an equivalence in $K$-theory  as seen in the proof of Proposition \ref{kgowpegrfwfwfref} and $K(-\otimes \Cl^{1})\simeq \Sigma K(-)$  we see that $I_{n+1}$ induces an equivalence in $K$-theory if and only if $I_{n}$ does so. In order to start the induction we consider the square 
 $$\xymatrix{  \bC_{(0)}(X, \hat H,\hat \phi) \ar[r]^{I_{(0)}}& C(X,\hat H,\hat \phi,0)  \\ \ar[u]\bC(X,H, \phi) \ar[r]^{I} &C(X, H, \phi)\ar[u]^{!} }$$
 The marked map can be identified with a left upper corner inclusion $C(X, H,\hat \phi)\to \Mat_{2}(C(X, H,\phi))$ with the standard grading. Again both vertical maps induce equivalences in $K$-theory. The map $I$ induces an equivalence in $K$-theory by the proof of   Theorem~\ref{fwefiwjfeiooi234234324434e}.
  \end{proof}

%
%
%
%

\begin{rem}\label{gig9egregre09greg}
In this remark we discuss the exterior tensor product.
We  consider two $G$-bornological coarse spaces $X$ and $X^{\prime}$. Furthermore, let $(H,\phi,n)$ and $(H^{\prime},\phi^{\prime},n^{\prime})$ be equivariant
$X$- and $X^{\prime}$-controlled  Hilbert spaces of degrees $n$ and $n^{\prime}$ which are determined on points.
We define the measure $\phi\hat\otimes \phi^{\prime}$ on $X\times X^{\prime}$ on 
$H\hat\otimes H^{\prime}$ by
$$(\phi\hat\otimes \phi^{\prime})(W):=\sum_{(x,x^{\prime})\in W} \phi(\{x\})\hat\otimes \phi^{\prime}(\{x^{\prime}\})\ .$$
 
Here we take advantage of the assumption that the factors are determined on points. In the general case, since we must define a measure on all subsets of $X\times X^{\prime}$, the construction  would be more complicated and would involve choices of suitable partitions of $X$ and $X^{\prime}$. 
 
 The graded Hilbert space $H\hat\otimes H^{\prime}$ carries an action of the graded algebra $\Cl^{n}\hat\otimes \Cl^{n^{\prime}}$ which is isomorphic to $\Cl^{n+n^{\prime}}$. Hence we get an equivariant  $(X\times X^{\prime})$-controlled Hilbert space $(H\hat\otimes H^{\prime},\phi\hat\otimes \phi^{\prime},n+n^{\prime})$ of degree $n+n^{\prime}$.
 
We furthermore get  a homomorphism of $C^{*}$-algebras 
$$C_{\lc}(X,H,\phi,n)\hat\otimes C_{\lc}(X^{\prime},H^{\prime},\phi^{\prime},n^{\prime})\to C_{\lc}(H\hat\otimes H^{\prime},\phi\hat\otimes \phi^{\prime},n+n^{\prime})$$
given by $A\hat\otimes A^{\prime}\mapsto A\hat\otimes A^{\prime}$
(here the domain is understood as an element of the abstract tensor product of $C^{*}$-algebras, and the image is an operator on $H\otimes H^{\prime}$). This map induces a homomorphism in $K$-theory groups \begin{equation}\label{g4pug45g4g}
\boxtimes: K_{\ell} (C_{\lc}(X,H,\phi,n))\otimes K_{\ell^{\prime}}(C_{\lc}(X^{\prime},H^{\prime},\phi^{\prime},n^{\prime}))\to K_{\ell+\ell^{\prime}}(C_{\lc}(H\hat\otimes H^{\prime},\phi\hat\otimes \phi^{\prime},n+n^{\prime}))
\end{equation}
which will be used in Section~\ref{goijrgoiregregreg}.
\hB
\end{rem}
\begin{rem}\label{erjogiopwerferfwrefwerfw}
The tensor product of $X$-controlled Hilbert spaces can  also be used to construct a lax symmetric monoidal refinement of the functor $\bC$ in \eqref{huihwiefewrfwfwrfwerfw}, and hence of $K\cX^{G}$ introduced in Definition \ref{girogegegergerg}, see  \cite[Sec. 4.2]{bunke2024localization}.
\hB
\end{rem}

%
%
%
%

\section{Proper complete \texorpdfstring{$\boldsymbol{G}$}{G}-manifolds are very proper}

 Let $M$ be a complete  Riemannian manifold with a proper and isometric  action of $G$. We consider $M$ as a $G$-bornological coarse space with the bornological and coarse structures induced by the Riemannian distance.
Recall Definition \ref{fpowjoefewfewf} of very properness.

\begin{prop}\label{kforererv}
$M$ is measurably very proper.
\end{prop}

\begin{proof}
We fix an open invariant entourage $V$ of $M$.
For a point $x$ in $M$ let $G_{x}$ denote the stabilizer subgroup of $x$.
There exists a linear action of $G_{x}$ on $\R^{n}$ and a $G_{x}$-invariant open neighbourhood $D_{x}$ of $0$ in $\R^{n}$ such that
a tubular neighbourhood   of the orbit $Gx$ in $M$ is equivariantly diffeomorphic to
$G\times_{G_{x}}D_{x}$. We let $U_{x}$ be the image of $\{1\}\times D_{x}$ under this diffeomorphism. We can assume that   the sets $U_{x}$ are $V$-bounded for all $x$ in $M$.

We consider the quotient $\bar M:=M/G$ (as a topological space) and write $\bar U_{x}$ for the image of $U_{x}$ under the natural projection $\pi:M\to \bar M$. The sets $(\bar U_{x})_{x\in M}$ form an open covering of $\bar M$.
We now use that $\bar M$ is $\sigma$-compact. 
Let $(\bar K_{n})_{n\in \nat}$ be an increasing exhaustion of $\bar M$ by compact subsets {with $\bar K_n \subseteq \inter(\bar K_{n+1})$ for all $n$.} For every integer $n$ we can choose a finite subset $\bar I _{n}$ of  
$\pi^{-1}(\bar K_{n}\setminus \inter(\bar K_{n-1}))$ such that
$\bar K_{n}\setminus \inter(\bar K_{n-1}) \subseteq \bigcup_{x\in \bar I _{n}} \bar U_{x}$.
For simplicity we can in addition assume that the $G$-orbits of the points in 
$\bar I _{n}$ are disjoint.
 For $x$ in $\bar I _{n}$ we define the $G_{x}$-invariant open subset $$V_{x}:=U_{x}\cap \pi^{-1}(\inter(\bar K_{n+1})\setminus \bar K_{n-2})\ .$$ 
We then define the subset 
$\bar I :=\bigcup_{n\in \nat} \bar I _{n}$ of $M$ and 
the $G$-closure
 $I :=G \bar I $. Furthermore, for  every $x$ in $\bar I$ and $g$ in $G$  define    the open subset
$V_{gx}:=gV_{x}$ (this is independent of the choice of $g$) of $M$.

We get an equivariant open covering $(V_{x})_{x\in I }$ such that
\begin{enumerate}
\item the family of subsets is uniformly $V$-bounded,
\item $I $ has finite stabilizers,
\item for every entourage $W$ of $M$ and $x$ in $I $ the set
$\{y\in I \:|\: W[V_{x}]\cap V_{y} \not=\emptyset\}$ is finite (indeed,
$\pi(W[V_{x}])$ is contained in $\bar K_{n}$ for some sufficiently large integer $n$),
\item for every $x$ in $I $ the point $x$ belongs to $V_{x}$ and is stabilized by $G_{x}$,
\item\label{greiogjoergergreg} for every bounded subset $B$ the set $\{y\in I \:|\:  V_{y} \cap B \not=\emptyset\}$ is finite (indeed, since we assume that $M$ is complete, the closure  of $B$ is compact and hence $\pi(B)$ is contained in $\bar K_{n}$ for some sufficiently large integer $n$), and
\item the set $I /G$ is countable.
\end{enumerate}
It remains to turn this covering into a partition.
 To this end we order the set
 $\bar I =I /G$ and choose an identification with $\nat$. For every integer $n$ we thus have a base point $x_{n}$ in the $n$'th orbit in $I $, namely the unique point of the orbit also belonging to $\bar I$. Using local finiteness of the covering, 
 for every  integer $n$   we can choose  a $G_{x_{n}}$-invariant neighbourhood $W_{x_{n}}$ of $x_{n}$ in $V_{x_{n}}$ which does not contain any other point of $I $. 
In a first step, for every integer $n$ we set 
$$V_{x_{n}}^{\prime}:=V_{x_{n}}\setminus \bigcup_{g\in G}\bigcup_{m\in \nat\setminus \{n\}} gW_{x_{m}}\ .$$ Then we extend this to an equivariant family of subsets 
$(V_{x}^{\prime})_{x\in J}$ by   setting
$V_{gx}^{\prime}:=gV_{x_{n}}^{\prime}$ for all integers $n$ and $g$ in $G$ (note that this is well-defined). Next we turn this family into a partition $(B_{x})_{x\in J}$ by setting
$$B_{gx_{n}}:=V_{gx_{n}}^{\prime}\setminus \bigcup_{m\in \nat, m< n } \bigcup_{h\in G} V_{hx_{m}}^{\prime}\ .$$
Since we have replaced the sets $V_{x}$ by the sets $V_{x}^{\prime}$ in the preceding step we have ensured that for every $x$ in $I$ we still have $x\in B_{x}$. 

We finally note that the sets $B_{x}$ for all $x$ in $I$ are measurable.
\end{proof}

Note that the proof of  Proposition \ref{kforererv} yields in  addition that we can prescribe the bound $V$ on the 
partition. This will be used below.

 Let $M$ be a complete  Riemannian manifold with a proper and isometric  action of $G$. 
We  consider $M$ as a $G$-bornological coarse space with the bornological and coarse structures induced by the Riemannian distance. In the following we show that  certain invariant open subsets $Z$ with the induced bornological coarse structures of $M$ are also very proper. The argument in \ref{greiogjoergergreg} above does not apply directly since bounded subsets of $Z$ need not be relatively compact. A simple idea is to intersect the partition constructed for $M$ with $Z$. But then in general the condition
 \ref{fpowjoefewfewf}.\ref{ojreoigjoergergerg} is violated. In the following we introduce an assumption  on $Z$ which ensures  \ref{fpowjoefewfewf}.\ref{ojreoigjoergergerg} in this procedure.

Recall  from \cite{equicoarse} that an invariant subset  $A$ of a $G$-bornological coarse space $X$ is called nice if  the  inclusion $A\to V[A]$  is an equivalence of $G$-bornological coarse spaces for every invariant entourage $V$ containing the diagonal.  The following notion introduces is a similar property. 
 
Let $V$ be an invariant entourage of $X$.

\begin{ddd}\label{rgklherogervfvvfgrege}\mbox{}
\begin{enumerate}
\item We say that $A$ is   nice for $V$  if for every
$x$ in $X$ and $a$ in $A\cap V[x]$ we have $G_{x}\subseteq G_{a}$.
\item We say that $A$ is a support if there  exists a cofinal set of invariant entourages $U$ of $X$
such that the $U$-thickening $U[A]$ of $A$ is nice for some entourage   (which may depend on $U$).
\end{enumerate}
\end{ddd}

Let $M$ be  as above and let $Z$ be an invariant open subset of $M$. We consider $Z$ as a $G$-bornological coarse space with the induced structures from $M$.

\begin{prop}\label{regioogergregreg}
If $Z$ is nice for some invariant open entourage $V$ of $M$, then
$Z$ is measurably very proper.
\end{prop}

\begin{proof}
  By Proposition \ref{kforererv} we can find a $V$-bounded equivariant measurable partition 
  $(B_{i})_{i\in I}$ of $M$ satisfying the conditions listed in Definition   \ref{fpowjoefewfewf}. We let $I^{\prime}:=\{i\in I\:|\: B_{i}\cap Z\not=\emptyset\}$ and
 define the equivariant partition
  $(Z\cap B_{i})_{i\in I^{\prime}}$ of $Z$.
  The only non-trivial condition to check is \ref{fpowjoefewfewf}.\ref{ojreoigjoergergerg}.
  By assumption, for every $i$ in $I^{\prime}$ there exists a point $x$ in $B_{i}$  such that $G_{i}\subseteq G_{x}$.
  Let $x^{\prime}$ be in $Z\cap B_{i}$. Then $x^{\prime}\in V[x]\cap Z
$. Since $Z$ is   nice for $V$ we have $G_{x}\subseteq G_{x^{\prime}}$. This implies
$G_{i}\subseteq G_{x^{\prime}}$ as required.
\end{proof}

\begin{ex}
 For $Z$ to be   nice for some open entourage $V$ of $M$ is not necessary for $Z$ being very proper.
We consider the action of the group $G:=\Z/3\Z$ by rotations   on $M:=\R^{2}$ with the standard Euclidean metric. We consider the open invariant subset $Z:=\R^{2}\setminus \{0\}$. If the unit ball $B(0,1)$ was
a member of our partition for $\R^{2}$ fixed by $G$, then  the problem is that $Z\cap B(0,1)$ does not contain any $G$-fixed point.
One can check that $Z$ is still very proper, but it is not   nice for any open   invariant entourage  of $\R^{2}$.
\hB
\end{ex}

\section{Construction of the index class}\label{goijrgoiregregreg}

Let $G$ be a discrete group and
 $M$ be a complete Riemannian manifold with a  proper   action of $G$ by isometries.  The Riemannian distance induces a $G$-bornological coarse structure on $M$. We use the symbol $M$ also  in order to denote the resulting $G$-bornological {coarse} space.

We consider a graded Hermitian vector bundle $E\to M$ with a right action of $\Cl^{n}$.
We form the $G$-Hilbert space $H_{0}:=L^{2}(M,E)$. It is graded and carries a right action by $\Cl^{n}$.
In order to turn this Hibert space into an equivariant $M$-controlled Hilbert space of degree $n$ we must construct an equivariant projection-valued measure $\phi_{0}$.  
 
By Proposition \ref{kforererv} the $G$-bornological coarse and $G$-measurable space $M$ is measurably very proper. 
We can thus choose an equivariant measurable partition of unity $(B_{i})_{i\in I}$ with the properties listed in Definition  \ref{fpowjoefewfewf}. In particular, by  Assumption \ref{fpowjoefewfewf}.\ref{ojreoigjoergergerg} we can choose an equivariant family $(b_{i})_{i\in I}$ of base points such that $b_{i}\in B_{i}$. We define
$$\phi_{0}:=\sum_{i\in I} \delta_{b_{i}} \chi(B_{i})\ ,$$
where $\chi(B_{i})$ is the multiplication operator on $H_{0}$  with the characteristic function of $B_{i}$, and $\delta_{b_{i}}$ is the Dirac measure at $b_{i}$.
We thus get an equivariant $M$-controlled Hilbert space $(H_{0},\phi_{0},n)$ of degree $n$ which is determined on points.

Note that the Hilbert space $H_{0}$ has a natural continuous control $\rho_{0}:C_{0}(M)\to B(H_{0})$,   where the elements of $C_{0}(M)$   act by multiplication operators. In this case $C_{\lc}(M,H_{0},\rho_{0},n)$ is the classical Roe algebra associated to $(H_{0},\rho_{0})$.

\begin{lem}
We have an equality
$C_{\lc}(M,H_{0},\rho_{0},n)=C_{\lc}(M,H_{0},\phi_{0},n)$.
\end{lem}
\begin{proof}
This uses that the members of the partition $(B_{i})_{i\in I}$ used to construct $\phi_{0}$ are uniformly bounded.
The argument for   the lemma is the same as the one \cite[Rem.~8.43]{buen} in the non-equivariant case. \end{proof}This comparison is relevant below  in order to apply the results of J. Roe \cite{Roe:2012fk}  which are formulated for the continuous controlled case.

We now consider an  invariant $\Cl^{n}$-linear Dirac operator $\Dirac $ acting on the sections of  $E$.
Such an operator is associated to a Dirac bundle structure on $E$ which  consists  in addition to the grading, the Hermitian metric, and the right $\Cl^{n}$-action, of a Clifford multiplication by tangent vectors, and of a connection satisfying some natural compatibility conditions. The Dirac operator $\Dirac $ is an essentially selfadjoint unbounded operator on $H_{0}$ with domain the smooth and compactly supported sections of $E$.  

We consider the $C^{*}$-algebra  $\hat C_{0}(\R)$ as in Section \ref{erkogpqwerfrewfrefrefw}.   For $\chi$ in $\hat C_{0}(\R)$ we define $\chi(\Dirac )$ in $B(H)$ using the functional calculus for essentially selfadjoint unbounded operators.

\begin{prop} 
The map
$\hat C_{0}(\R)\to B(H)$ given by $\chi\mapsto \chi(\Dirac )$
is a   homomorphism of graded $C^{*}$-algebras $I(\Dirac ):  \hat C_{0}(\R)\to C_{\lc}(M,H_{0},\phi_{0},n)$.
\end{prop}
\begin{proof}
This is a basic and well-known fact  in coarse index theory.
If $\chi$ in $\hat C_{0}(\R)$ is  such that the Fourier transform
$\hat \chi$ is smooth and has compact support, then
$\chi(\Dirac)$ is $G$-invariant, locally compact and controlled. One can approximate the elements of $\hat C_{0}(\R)$ uniformly by such functions.
The homomorphism $I(\Dirac)$ preserves the grading since $\Dirac$ is odd.
\end{proof}

\begin{ex}\label{erigjorewgwrferfrwefwref}
We consider the $\Cl^{1}$-linear Dirac operator  on $\R^{1}$ acting on $L^{2}(\R)\hat \otimes \Cl^{1}$ by 
$\Dirac:=\sigma \partial_{x}  $, where $\sigma$ is the odd  anti-selfadjoint generator of   $\Cl^{1}$ with $\sigma^{2}=-1$ acting on $\Cl^{1}$ by left multiplication. The operator $\Dirac$ commutes with the action of $\Cl^{1}$ on itself by right multiplication. We get a homomorphism $I(\Dirac):\hat C_{0}(\R)\to C_{lc}(\R^{1},L^{2}(\R),\rho_{0})\hat \otimes \Cl^{1}$ which represents a class $\beta$  in
$ K_{0}(C_{lc}(\R^{1},L^{2}(\R),\rho_{0})\hat \otimes \Cl^{1})\simeq  K_{-1}(\Cl^{1})$.
The equivalence $K(A)\stackrel{\simeq}{\to} \Sigma K(A\hat \otimes \Cl^{1})$ in Example \eqref {jiefgoerfqewefqewd} is induced by multiplication with $\beta$.
\hB
 
\end{ex}

 Let $B$ be a subset of $M$. The Weizenboeck formula 
 $$\Dirac^{2}=\nabla^{*}\nabla+R$$ determines a selfadjoint bundle endomorphism $R$ in $C^{\infty}(M,\End(E)^{sa})$, where $\nabla$ denotes the connection on the Dirac bundle $E$.
 
\begin{ddd}\label{goiergregg}
$\Dirac $ is uniformly locally positive on $B$ if there exists  a positive real number $c$  such that for every $x$ in $B$ we have
 $c \cdot \id_{E_{x}}\le R(x)$.
\end{ddd}
We call $c$ a lower bound for $\Dirac$ on $B$.

By \cite[Def. 3.5]{equicoarse} a (equivariant)  big family $\cY=(Y_{i})_{i\in I}$ in a $G$-bornological coarse space $X$ is a family of invariant subsets of $X$ indexed by a poset such that $i\le j$ implies $Y_{i}\subseteq Y_{j}$, and such that 
for every coarse entourage $U$ of $X$ and $i$ in $I$ there exists $j$ in $I$ with $U[Y_{i}]\subseteq Y_{j}$, where $U[Y_{i}]$ denotes the $U$-thickening of $Y_{i}$.

We define the Roe algebra $C_{\lc}(\cY,H_{0},\phi_{0},n)$ of a big family $\cY$ in $M$ as the $C^{*}$-subalgebra of the Roe algebra $C_{\lc}(M,H_{0},\phi_{0},n)$ generated by operators which are supported on members of the family $\cY$. Recall that $\{A\}$ denotes the big family of thickenings of the invariant subset $A$ of $M$. {In the ungraded case} we have the equality $C_{\lc}(\{A\},H_{0},\phi_{0})= C^{*}(A\subseteq M)$, where the right-hand side is Roe's notation \cite{Roe:2012fk}.

 We now assume that $A$ is a $G$-invariant subset of $M$ such that $\Dirac $ is 
uniformly locally positive on the complement $M\setminus A$ with lower bound $c^2$.

Let $\chi$ be in $C_{0}(\R)$.
\begin{prop}
\label{fiowefewfewf}
 If   $\supp(\chi)\in (-c,c)$, then
 we have 
 $\chi(\Dirac )\in C_{\lc}(\{A\},H_{0},\phi_{0},n)$.
\end{prop}
\begin{proof}
This has been shown by J. Roe \cite[Lem.~2.3]{Roe:2012fk}.
\end{proof}


By Proposition \ref{fiowefewfewf} the homomorphism $I(\Dirac)$ restricts to a homomorphism
$$I_{c}(\Dirac ):\hat C_{0}((-c,c))\to C_{\lc}(\{A\},H_{0},\phi_{0},n)\ .$$

This homomorphism represents the primary index class (cf. Zeidler \cite[Def. 4.1]{MR3551834})
\begin{equation}\label{g4poi450pg45g45g}
i(\Dirac,\on(A)) \:\:\mbox{in} \:\: K_{0}(C_{\lc}(\{A\},H_{0},\phi_{0},n))
\end{equation}
in the $K$-theory (for graded $C^{*}$-algebras). It is independent of the choice of $c$.

 Let $(H ,\phi )$ be an ample equivariant $M$-controlled Hilbert space. By Lemma \ref{okgpergergreg} we can choose
 an equivariant controlled and $\Cl^{n}$-linear unitary embedding $$U:(H_{0},\phi_{0},n)\to (\hat H ,\hat \phi ,n)\ .$$ We can assume that this embedding has an ample complement. To this end, if necessary,  we  compose it further with the first summand embedding $(\hat H ,\hat \phi ,n)\to (\hat H\oplus \hat H ,\hat \phi\oplus \hat \phi ,n)$ and change notation afterwards.
  By Lemma  \ref{irjfiorregregregre} any two such embeddings are homotopic. 

The embedding $U$ induces a homomorphism of $C^{*}$-algebras
 $$C_{\lc}({\{A\}},H_{0},\phi_{0},n)\to C_{\lc} ({\{A\}},\hat H  ,\hat \phi ,n) \ .$$ 
 This homomorphism of $C^{*}$-algebras induces a homomorphism of $K$-theory spectra
 $$K(C_{\lc}({\{A\}},H_{0},\phi_{0},n))\to K(C_{\lc} ({\{A\}},\hat H  ,\hat \phi ,n))$$
 which by homotopy invariance of $K$  is independent of the choice of $U$.
 We now assume that $A$ is a support (Definition \ref{rgklherogervfvvfgrege}).
Then
 we have  an equivalence 
 $$\kappa:K(C_{\lc}({\{A\}},H_{0},\phi_{0},n)) \stackrel{\simeq}{\leftarrow} K(C({\{A\}},H_{0},\phi_{0},n))\stackrel{\simeq}{\to}  \KX^{G}_{(n)}({\{A\}})\ ,$$
 where the first is 
given by  Proposition \ref{lem:sdfbi23} (applied to the members $U[A]$ of $\{A\}$ for a   cofinal subfamily of invariant entourages $U$ of  $X$ such that   $U[A]$ is nice for some invariant entourage and hence very proper by Proposition \ref{regioogergregreg}), and the second  is  Theorem \ref{fwefiwjfeiooi234234324434en}.    

%
\begin{ddd}[Coarse index class with support]\label{gioregegergregr}
The index class of $\Dirac $ with support on $A$ is defined by $$\Ind(\Dirac,\on(A)) :=\kappa(i(\Dirac,\on(A)))\:\:\mbox{in}\:\: \pi_{0}\KX^{G}_{(n)}({\{A\}})\ .$$
\end{ddd}

\begin{rem}
In view of the equivalence 
$\KX^{G}_{(n)} \simeq \Omega^{n} K\cX^{G} $ from Proposition \ref{kgowpegrfwfwfref}
we can consider
the index class also as an element 
$\Ind(\Dirac,\on(A))$ in $K\cX_{n}^{G}({\{A\}})$. \hB
\end{rem}

Let us discuss now the compatibility of the index classes of Dirac operators with products. 
Let  $\Dirac$ be an invariant  Dirac operator of degree $n$ on a complete Riemannian manifold $M$ with isometric proper $G$-action {and such that $\Dirac$} is uniformly locally positive  outside of an invariant subset $A$ of $M$ which is a support. Let $\Dirac^{\prime}$, $M^{\prime}$, $A^{\prime}$, $n^{\prime}$ be similar  data.  
 Then we have classes
 $$i(\Dirac,\on(A))\:\:\mbox{in}\:\: K_{0}(C_{lc}({\{A\}},H_{0},\phi_{0},n))$$
and 
 $$i(\Dirac^{\prime},\on(A^{\prime}))\:\:\mbox{in}\:\: K_{0}(C_{lc}({\{A^{\prime}\}} ,H^{\prime}_{0},\phi^{\prime}_{0},n^{\prime}))\ .$$
We can form the invariant  Dirac operator 
 $\Dirac\hat \otimes 1+1\hat\otimes \Dirac^{\prime}$ on $M\times M^{\prime}$ of degree $n+n^{\prime}$. It is  uniformly locally positive  outside
 $A\times A^{\prime}$ which is again a support.
 We therefore get a class 
 $$i(\Dirac\hat\otimes 1+1\hat\otimes \Dirac^{\prime},\on(A\times A^{\prime}))\:\:\mbox{in}\:\: K_{0}(C_{\lc}({\{A\times A^{\prime}\}},H_{0}^{\prime\prime} ,\phi^{\prime\prime}_{0},n+n^{\prime}))\ .$$
Recall that $H_{0}=L^{2}(M,E)$, $H_{0}^{\prime}=L^{2}(M^{\prime},E^{\prime})$ and $H_{0}^{\prime\prime}=L^{2}(M\times M^{\prime},E\otimes E^{\prime})$. There is a canonical isomorphism of Hilbert spaces
$$u:H_{0}\hat \otimes H_{0}^{\prime}\stackrel{\cong}{\to} H_{0}^{\prime\prime}\ .$$ 
It preserves the grading, is $\Cl^{n+n^{\prime}}$-linear, and it is controlled if we equip the domain with  the measure $\phi_{0}\hat\otimes \phi_{0}^{\prime}$ and the target with the measure $\phi_{0}^{\prime\prime}$ (see Remark \ref{gig9egregre09greg}), where $\phi_{0}$, $\phi_{0}^{\prime}$ and   $\phi_{0}^{\prime\prime}$ are  constructed by the procedure explained above.
In particular, we have an isomorphism of graded $C^{*}$-algebras
$$u_{*}: C_{\lc}(\{A\times A^{\prime}\} ,H_{0}\hat\otimes H_{0}^{\prime},\phi_{0}\hat\otimes \phi^{\prime}_{0},n+n^{\prime})\to C_{\lc}(\{A\times A^{\prime}\},H_{0}^{\prime\prime},\phi_{0}^{\prime\prime},n+n^{\prime})\ , \quad B\mapsto uBu^{*}\ .$$
Recall the exterior product \eqref{g4pug45g4g}.  
\begin{prop}\label{gropgergergerge}
We have 
$$u_{*}\big( i(\Dirac,\on(A))\boxtimes i(\Dirac^{\prime},\on(A^{\prime})) \big) = i(\Dirac\hat \otimes 1+1\hat\otimes \Dirac^{\prime},\on(A\times A^{\prime}))\ .$$
\end{prop}

\begin{proof}
The proof is the same as the proof of Zeidler \cite[Thm. 4.14]{MR3551834}. Zeidler assumes free actions, but this not relevant for the argument. 
\end{proof}

\begin{rem}
If we interpret the index classes   in $K\cX^{G}$, then in terms of the lax symmetric monoidal structure 
$$K\cX^{G}(M)\otimes K\cX^{G}(M')\to K\cX^{G}(M\times M')$$ of $K\cX^{G}$ mentioned in Remark \ref{erjogiopwerferfwrefwerfw}
  the Proposition \ref{gropgergergerge} says that
$ i(\Dirac\hat\otimes 1+1\hat\otimes \Dirac^{\prime},\on(A\times A^{\prime}))$ in $K\cX^{G}_{n+n'}(M\times M')$ is the image of the product
$ i(\Dirac,\on(A))\otimes i(\Dirac^{\prime},\on(A^{\prime})$ in $K\cX^{G}_{n}(M)\otimes K\cX^{G}_{n'}(M)$.
\hB
\end{rem}

\section{A coarse relative index theorem}\label{groigeorgegerg}

In this section we prove a relative index theorem for the coarse index with support. It is the technically precise expression that two Dirac operators on two manifolds which are isomorphic to each other on some subset have the same index if the index is considered as an element in the relative coarse $K$-homology of the manifold relative to the complements of the respective subsets. In order to compare the indices we use the   excision morphism  between the   relative coarse $K$-homology groups.
We will actually show a slightly more general theorem which takes into account that an equivariant isometry between two subsets is not necessarily an isomorphism of $G$-bornological coarse spaces if their structures are induced from the respective ambient manifolds.

 We start with the description of the data needed to state the coarse relative index theorem.  
Let $G$ be a discrete group.  We consider a 
  complete Riemannian manifold $M$ with a proper  action of $G$ by isometries and  with an invariant Dirac operator $\Dirac_{M}$ of degree $n$. We assume that there is an   invariant subset 
  $A$ of $M$ which is a support and such that $\Dirac_{M}$ is uniformly locally positive on the complement of $A$.   
   
Let $Z$ be an open and very proper (e.g., nice for some invariant open entourage of $M$, see Proposition \ref{regioogergregreg}) invariant  subset of $M$ and let $Z^{c}$ denote its complement. 
{Note that since $A$ is a support, there is a cofinal set of invariant entourages $U$ such that $Z\cap U[A]$ is measurably very proper (cf.~Proposition~\ref{regioogergregreg}).}
We  have the big family $\{Z^{c}\}\cap Z$ in $Z$.


%

Let $M^{\prime}$, $A^{\prime}$, $\Dirac_{M^{\prime}}$ and $Z^{\prime}$ be similar data.

We assume that there is an equivariant diffeomorphism  $i:Z\stackrel{\cong}{\to} Z^{\prime}$ which preserves the Riemannian metric.  
We assume that 
$i$ also induces a morphism of $G$-bornological coarse spaces where the bornological coarse structures on $Z$ or $Z^{\prime}$ are induced from $M$ or $M^{\prime}$, respectively.

\begin{rem}
We note that the coarse structure on $Z$ does not only depend on the metric of $Z$. Two points which are very distant to each other in $Z$ might be actually close to each other   in $M$ since they might be connected by a  short path leaving $Z$.
So our assumption is that for points $x,y$ in $Z$ the  distance between the image  points $i(x)$ and $i(y)$ in $Z^{\prime}$ measured in the geometry of $M^{\prime}$ can be bounded in terms of  the distance between $x$ and $y $ measured in $M$.\hB
\end{rem}

In order to state our compatibility assumptions for the big families we use the following language.
 We consider sets $X$ and $X^{\prime}$ and a map of sets $f:X\to X^{\prime}$.
 Let $\cY=(Y_{i})_{i\in I}$ and $\cY^{\prime}=(Y^{\prime}_{i^{\prime}})_{i^{\prime}\in I^{\prime}}$ be filtered families of subsets  of  $X$ and $X^{\prime}$, respectively. We say that
 $f$ induces a morphism $\cY\to \cY^{\prime}$ if for every $i$ in $I$ there exists $i^{\prime}$ in $I^{\prime}$ such that
 $f(Y_{i})\subseteq Y^{\prime}_{i^{\prime}}$.
 
We now continue with the assumptions for the coarse relative index theorem.
   We  assume that the map $i:Z\to Z^{\prime}$
   induces  a morphism  between the big families $\{A\}\cap Z$ and  $\{A^{\prime}\}\cap Z^{\prime}$, and a morphism  between the big families
$\{Z^{c}\}\cap Z$ and $\{Z^{\prime,c}\}\cap Z^{\prime}$.

\begin{rem}
The family big family $\{Z^{c}\}\cap Z$ is the coarse geometric version of the boundary of $Z$.
\hB
\end{rem}

 We finally  assume that the isometry $i$ is
 covered by an equivariant isomorphism of Dirac operators 
$(\Dirac_{M})_{|Z}\cong (\Dirac_{M^{\prime}})_{|Z^{\prime}}$.

\begin{ddd}
We call  data  satisfying the {above} assumptions a coarse relative index situation.
\end{ddd}

We have a diagram of horizontal pair-fibre sequences of spectra \begin{align}
\mathclap{\xymatrix{
\KX^{G}_{(n)}(\{Z^{c}\}\cap \{A\})\ar[r]&\KX_{(n)}^{G}( \{A\})\ar[r]&\KX_{(n)}^{G}(  \{A\},\{Z^{c}\}\cap \{A\})\\
\KX_{(n)}^{G}(Z\cap \{Z^{c}\}\cap \{A\})\ar[d]\ar[r]\ar[u]&\KX_{(n)}^{G}(Z\cap \{A\})\ar[r]\ar[u]\ar[d]& \KX_{(n)}^{G}(Z\cap \{A\},Z\cap \{Z^{c}\}\cap \{A\})\ar[u]_{\simeq}\ar[d]\\
\KX_{(n)}^{G}(Z^{\prime}\cap \{Z^{\prime,c}\}\cap \{A^{\prime}\}) \ar[d]\ar[r]&\KX_{(n)}^{G}(Z^{\prime}\cap \{A^{\prime}\})\ar[r]\ \ar[d]& \KX_{(n)}^{G}(Z^{\prime}\cap \{A^{\prime}\} ,Z^{\prime}\cap \{Z^{\prime,c}\}\cap \{A^{\prime}\})\ar[d]^{\simeq}\\
\KX_{(n)}^{G}(\{Z^{\prime,c}\}\cap \{A^{\prime}\})\ar[r]&\KX_{(n)}^{G}(  \{A^{\prime}\})\ar[r]&\KX_{(n)}^{G}(  \{A^{\prime}\},\{Z^{\prime,c}\}\cap \{A^{\prime}\})
}}\notag\\
\mbox{}\label{hgiregooerg8989}
\end{align}
The lower and the upper left square are push-out squares by excision. This explains the lower and upper right vertical equivalences.
The middle vertical  morphisms  are induced by the morphism $i:Z \to Z^{\prime}$.

\begin{ddd}
The morphism induced by the right column in \eqref{hgiregooerg8989}
$$e:\KX_{(n)}^{G}(\{A\} ,\{Z^{c}\}\cap \{A\}) \to \KX_{(n)}^{G}( \{A^{\prime}\},\{Z^{\prime,c}\}\cap \{A^{\prime}\})$$
is called the excision morphism associated to the coarse relative index situation.
\end{ddd}

We let
$$\overline{\Ind(\Dirac_{M},\on(A))}\:\:\mbox{in}\:\: \pi_{0}\KX^{G}_{(n)}(\{A\},\{Z^{c}\}\cap \{A\})$$   denote the image of  $ \Ind(\Dirac_{M},\on(A))$ under the natural map $$\KX_{(n)}^{G} (\{A\} )\to \KX_{(n)}^{G} (\{A\},\{Z^{c}\}\cap \{A\}) \ .$$ 
We similarly define the class
$$\overline{\Ind(\Dirac_{M^{\prime}},\on(A^{\prime}))}\:\:\mbox{in}\:\:\pi_{0}\KX_{(n)}^{G}(\{A^{\prime}\},\{Z^{\prime,c}\}\cap \{A^{\prime}\})\ .$$

\begin{theorem}[Coarse relative index theorem]\label{roigeorgergergerg}
We have the equality $$e(\overline{\Ind(\Dirac_{M},\on(A))})=\overline{\Ind(\Dirac_{M^{\prime}},\on(A^{\prime}))}\ .$$
\end{theorem}
\begin{proof}
The proof has two parts. This first is the comparison
of the index classes of the Dirac operators as $K$-theory classes  of respective quotients of Roe algebras. This part is probably well-known to the experts in the field.

The second part is the transition from the $K$-theories of  Roe algebras to the $K$-theory of Roe categories which are in the background of the construction of the equivariant coarse $K$-homology functor $\KX^{G}$.

We start with the basic analytic facts.
We consider the Hilbert spaces $H_{0}:=L^{2}(M,E)$
and $H_{0}^{\prime}:=L^{2}(M^{\prime},E^{\prime})$, where $E$ and $E^{\prime}$ denote the corresponding Dirac bundles underlying $\Dirac_{M}$ and $\Dirac_{M^{\prime}}$. Also, we let $\phi_{0}$ and $\phi^{\prime}_{0}$ be the projection-valued measure defined by multiplication with the characteristic functions of Borel measurable subsets. Note that the pairs $(H_{0},\phi_{0})$ and $(H_{0}^{\prime},\phi_{0}^{\prime})$ are not equivariant $X$-controlled Hilbert spaces   since the measures are only defined on Borel subsets.
But this measurable control suffices to define the Roe algebras appearing in the proof. We abuse notation and use the same symbols as in the case of equivariant $X$-controlled Hilbert spaces.
All occuring subsets of $M$ or $M^{\prime}$ in the following are assumed to be measurable.

The  {isomorphism between the objects  over $Z$ and $Z^{\prime}$, respectively,  given by the data of a relative index situation},  induces an  isometry $u:H_{0}(Z)\to H^{\prime}_{0}(Z^{\prime})$ such that $i_{*}\phi_{0|Z}=u^{*}\phi_{0|Z^{\prime}}^{\prime} u$.  


It is known from proofs of excision for coarse $K$-homology (e.g., \cite[Proof of Prop.~8.82]{buen}, \cite[Lem. 6.17]{bu}) that the inclusion
$$C_{\lc}(Z\cap \{A\},H_{0},\phi_{0},n)\to C_{\lc}(\{A\},H_{0},\phi_{0},n)$$
induces an isomorphism of $C^{*}$-algebras
\begin{equation}\label{cweihoi5gg4g}
\frac{C_{\lc}(Z\cap \{A\},H_{0},\phi_{0},n)}{C_{\lc}(Z\cap \{Z^{c}\}\cap\{A\},H_{0},\phi_{0},n)} \xrightarrow{\cong}
\frac{C_{\lc}(\{A\},H_{0},\phi_{0},n)}{C_{\lc}(\{Z^{c}\}\cap\{A\},H_{0},\phi_{0},n)}\ .
\end{equation}
Here we define the Roe algebra  of a big family as the $C^{*}$-subalgebra of the Roe algebra of the ambient space generated by operators which are supported on members of the family.
The inverse of the isomorphism \eqref{cweihoi5gg4g} is given by 
 $[F]\mapsto  [\phi_{0}(Z)F\phi_{0}(Z)]$. We have a similar isomorphism for the primed objects.
 The isometry $u$ induces a morphism of $C^{*}$-algebras
 $$[u]: \frac{C_{\lc}(Z\cap \{A\},H_{0},\phi_{0},n)}{C_{\lc}(Z\cap \{Z^{c}\}\cap\{A\},H_{0},\phi_{0},n)}\to \frac{C_{\lc}(Z^{\prime}\cap \{A^{\prime}\},H_{0}^{\prime},\phi_{0}^{\prime},n)}{C_{\lc}(Z^{\prime}\cap \{Z^{\prime,c}\}\cap\{A^{\prime}\},H_{0}^{\prime},\phi_{0}^{\prime},n)} $$
 given by
 $[F]\mapsto [uFu^{*}]$.
 In general it is not an isomorphism, because the control conditions in the target might be weaker than in the source. 
  
  Let $c$ in $(0,\infty)$ be a lower bound for the uniform local positivity of $\Dirac_{M}$ and $\Dirac_{M^{\prime}}$ on the complements of the sets $A$ and $A^{\prime}$, respectively, see Definition \ref{goiergregg}.
  We consider a function  $\chi$ in $C_{0}^{\infty}((-c,c))$.  
\begin{lem}\label{fewoijfowefwefwef}
We have an equality
 $$[u]([\phi_{0}(Z)\chi(\Dirac_{M})\phi_{0}(Z)])=[\phi^{\prime}_{0}(Z^{\prime})\chi(\Dirac_{M^{\prime}})\phi^{\prime}_{0}(Z^{\prime})]\ .$$
\end{lem}
 
\begin{proof}
Let $$\Delta:=u\phi_{0}(Z)\chi(\Dirac_{M})\phi_{0}(Z)u^{*}-\phi^{\prime}_{0}(Z^{\prime})\chi(\Dirac_{M^{\prime}})\phi^{\prime}_{0}(Z^{\prime})\ .$$
We must show that $\Delta\in C_{\lc}(Z^{\prime}\cap \{Z^{\prime,c}\}\cap\{A^{\prime}\},H_{0}^{\prime},\phi_{0}^{\prime},n)$.

Note that
$\Delta\in C_{\lc}(Z^{\prime}\cap \{A^{\prime}\},H_{0}^{\prime},\phi_{0}^{\prime},n)$ by Proposition \ref{fiowefewfewf}.
 We fix   $\epsilon$ in $(0,\infty)$. Then we fix $R$ in $(0,\infty)$ so large such that
$$\frac{1}{\sqrt{2\pi}}\int_{\R\setminus [-R,R]} |\hat \chi(t)| \, dt\le\frac{ \epsilon}{4}\ .$$ 
We have 
$U^{\prime}_{R}[Z^{\prime}\setminus U_{R}^{\prime}[Z^{\prime,c}]]\cap Z^{\prime,c}=\emptyset$, where $U_{R}^{\prime}:=\{(x^{\prime},y^{\prime})\in M^{\prime}\times M^{\prime}\:|\: \dist_{M^{\prime}}(x^{\prime},y^{\prime})\le R \}$. We use the 
  equality
$$\chi(\Dirac_{M})=\frac{1}{\sqrt{2\pi}}\int_{\R} e^{it\Dirac_{M}} \hat \chi(t) \, dt\ .$$
We get 
  $$\phi^{\prime}_{0}(Z^{\prime}\setminus U_{R}^{\prime}[Z^{\prime,c}])\Delta=\frac{1}{\sqrt{2\pi}}\int_{\R\setminus [-R,R]} \phi^{\prime}_{0}(Z^{\prime}\setminus U_{R}^{\prime}[Z^{\prime,c}]) (ue^{it\Dirac_{M}} u^{*}-e^{it\Dirac_{M^{\prime}}}) \phi_{0}^{\prime}(Z^{\prime}) \hat \chi(t) \, dt\ ,$$
where we can omit the interval $[-R,R]$ from the domain of integration since the integrand vanishes there identically {firstly} by the unit propagation speed of the wave groups $t\mapsto e^{it\Dirac_{M}}$ and $t\mapsto e^{it\Dirac_{M^{\prime}}}$, {and secondly because by assumption we have $i_{*}\phi_{0|Z}=u^{*}\phi_{0|Z^{\prime}}^{\prime} u$ and $i$ is covered by an equivariant isomorphism of Dirac operators $(\Dirac_{M})_{|Z}\cong (\Dirac_{M^{\prime}})_{|Z^{\prime}}$.}
 In view of the choice of $R$ and the unitarity of the wave operators this implies the first one of the following two estimates:
\[\|\phi^{\prime}_{0}(Z^{\prime}\setminus U_{R}^{\prime}[Z^{\prime,c}])\Delta   \| \le 
  \frac{\epsilon}{2}\ ,\quad \|\Delta \phi^{\prime}_{0}(Z^{\prime}\setminus U_{R}^{\prime}[Z^{\prime,c}])\| \le 
  \frac{\epsilon}{2} \ .\]
The second inequality follows by considering adjoints with $\chi$ replaced by {the function} given by $t\mapsto \bar \chi(-t)$.

We now write
\begin{align*}
\mathclap{
\Delta -\phi_{0}^{\prime}( Z^{\prime}\cap  U_{R}^{\prime}[Z^{\prime,c}])\Delta \phi^{\prime}_{0}(  Z^{\prime}\cap  U_{R}^{\prime}[Z^{\prime,c}])=  \Delta \phi^{\prime}_{0}(Z^{\prime}\setminus  U_{R}^{\prime}[Z^{\prime,c}])+  \phi^{\prime}_{0}(  Z^{\prime}\setminus U_{R}^{\prime}[Z^{\prime,c}]) \Delta \phi^{\prime}_{0}(  Z^{\prime}\cap V^{\prime}[Z^{\prime,c}])
}
\end{align*}
and conclude that
$$\| \Delta- \phi_{0}^{\prime}( Z^{\prime}\cap  U_{R}^{\prime}[Z^{\prime,c}])\Delta \phi^{\prime}_{0}(  Z^{\prime}\cap  U_{R}^{\prime}[Z^{\prime,c}])\|\le \epsilon\ .$$
But
$\phi_{0}^{\prime}( Z^{\prime}\cap  U_{R}^{\prime}[Z^{\prime,c}])\Delta \phi^{\prime}_{0}(  Z^{\prime}\cap  U_{R}^{\prime}[Z^{\prime,c}])$ belongs to 
$C_{\lc}(Z^{\prime}\cap \{Z^{\prime,c}\}\cap\{A^{\prime}\},H_{0}^{\prime},\phi_{0}^{\prime},n)$.
Since $\epsilon$ can be taken arbitrarily small
we conclude that $\Delta\in C_{\lc}(Z^{\prime}\cap \{Z^{\prime,c}\}\cap\{A^{\prime}\},H_{0}^{\prime},\phi_{0}^{\prime},n)$.
\end{proof}
 
Let us continue with the proof of Theorem~\ref{roigeorgergergerg}. We let $\overline{ i(\Dirac_{M},\on(A))}$  be the image  of the index class
\eqref{g4poi450pg45g45g} under the {composition} 
\begin{align*}
K_{0}(C_{\lc}(\{A\},H_{0},\phi_{0},n)) & \to K_{0} \Big(\frac{C_{\lc}(\{A\},H_{0},\phi_{0},n)}{C_{\lc}(\{Z^{c}\}\cap \{A\},H_{0},\phi_{0},n)}\Big)\\
& \xrightarrow{\cong} K_{0}\Big(\frac{C_{\lc}(Z\cap \{A\},H_{0},\phi_{0},n) }{C_{\lc}(Z\cap \{Z^{c}\}\cap \{A\},H_{0},\phi_{0},n)}\Big)\ .
\end{align*}
We define
$\overline{i(\Dirac_{M^{\prime}},\on(A^{\prime}))}$ similarly.
Then  Lemma \ref{fewoijfowefwefwef} shows that
\begin{equation}\label{knferlknlefnlrefref34r}
[u]_{*}(\overline{ i(\Dirac_{M},\on(A ))})=\overline{ i(\Dirac_{M^{\prime}},\on(A^{\prime}))}\ .
\end{equation}
The rest of the proof of the {Theorem~\ref{roigeorgergergerg}} consists of the (tedious and not very exiting) verification that this equality  implies the equality claimed in the theorem.
To this end we must compare the morphism 
$[u]_{*}$  with the excision morphism.

In a first step we choose ample $X$-controlled Hilbert spaces  $(H,\phi)$ and $(H^{\prime},\phi^{\prime})$ on $M$ and on $M^{\prime}$. 
We can assume that $(H(Z),\phi_{|Z})$ is ample on $Z$, too.

We can and will furthermore assume that $H(Z)=H^{\prime}(Z^{\prime})$. To this end we can replace the original spaces $(H,\phi)$ and $(H^{\prime},\phi^{\prime})$ by
$(H\oplus H^{\prime}(Z^{\prime}),\phi\oplus i^{-1}_{*}\phi^{\prime}_{|Z^{\prime}})$ and $(H^{\prime}\oplus H(Z),\phi^{\prime}\oplus i_{*}\phi_{|Z})$, respectively, and
then identify the corresponding subspaces.

Using Lemma \ref{okgpergergreg} we choose equivariant controlled  and $\Cl^{n}$-linear embeddings with ample complements
$v:H_{0}\to \hat H$ and $v^{\prime}:H_{0}^{\prime}\to \hat H^{\prime}$ such that
$v^{\prime}\circ u=v_{|H(Z)}$. This situation can be ensured by first defining $v$, then defining 
$v^{\prime}$ by this formula on $H^{\prime}(Z^{\prime})$ and finally extending it to all of $H_{0}^{\prime}$. 

We then have a commutative diagram of $C^{*}$-algebras
\[\mathlarger{\xymatrix{
\frac{C_{\lc}(\{A\},H_{0},\phi_{0},n)}{C_{\lc}(\{Z^{c}\}\cap\{A\},H_{0},\phi_{0},n)}\ar[r]&\frac{C(\{A\},\hat H,\hat \phi,n)}{C(\{Z^{c}\}\cap\{A\},\hat H,\hat \phi,n)}\\
\frac{C_{\lc}(Z\cap \{A\},H_{0},\phi_{0},n)}{C_{\lc}(Z\cap \{Z^{c}\}\cap\{A\},H_{0},\phi_{0},n)}\ar[r]^{[v]}\ar[u]_{\cong}\ar[d]_{[u]} &\frac{C(Z\cap \{A\},\hat H,\hat \phi,n)}{C(Z\cap \{Z^{c}\}\cap\{A\},\hat H,\hat \phi,n)}\ar[d]_{f}\ar[u]_{\cong}\\
\frac{C_{\lc}(Z^{\prime}\cap \{A^{\prime}\},H^{\prime}_{0},\phi^{\prime}_{0},n)}{C_{\lc}(Z^{\prime}\cap \{Z^{\prime,c}\}\cap\{A^{\prime}\},H^{\prime}_{0},\phi^{\prime}_{0},n)} \ar[r]^{[v^{\prime}]}\ar[d]^{\cong}&\frac{C(Z^{\prime}\cap \{A^{\prime}\},\hat H^{\prime},\hat \phi^{\prime},n)}{C(Z^{\prime}\cap \{Z^{\prime,c}\}\cap\{A^{\prime}\},\hat H^{\prime},\hat \phi^{\prime},n)}\ar[d]^{\cong}\\
\frac{C_{\lc}(\{A^{\prime}\},H^{\prime}_{0},\phi^{\prime}_{0},n)}{C_{\lc}(\{Z^{\prime,c}\}\cap\{A^{\prime}\},H^{\prime}_{0},\phi^{\prime}_{0},n)}\ar[r]&\frac{C(\{A^{\prime}\},\hat H^{\prime},\hat \phi^{\prime},n)}{C(\{Z^{\prime,c}\}\cap\{A^{\prime}\},\hat H^{\prime},\hat \phi^{\prime},n)}
}}\]
where the horizontal homomorphisms are induced by the embeddings $v$ and $v^{\prime}$, respectively. We further use the very properness of $M$  and of $Z$ (Proposition \ref{kforererv} for $M$ and  an assumption for  $Z$)    in order to drop the subscript $-_{\lc}$ on the right column by Proposition \ref{lem:sdfbi23}. 
The morphism $f$ is induced by the inclusion $C(Z\cap \{A\},\hat H,\hat \phi,n)\to C(Z^{\prime}\cap \{A^{\prime}\},\hat H^{\prime},\hat \phi^{\prime},n)$.

By \eqref{knferlknlefnlrefref34r} we have 
\begin{equation}\label{knferlknlefnlrefref34r1}
 [v^{\prime}]_{*}( \overline{ i(\Dirac_{M^{\prime}}, \on (A^{\prime} )}))=f_{*}([v]_{*}(\overline{ i(\Dirac_{M }, \on (A ))}))\ .
\end{equation}


We now have a commutative diagram {(see the proof  of Theorem ~\ref{fwefiwjfeiooi234234324434en} for notation)}
\begin{align}\label{ferfpokp43rgr}
\mathlarger{\mathclap{\xymatrix{
K\big(\frac{C(\{A\},H,\phi,n)}{C(\{Z^{c}\}\cap\{A\},H,\phi,n)}\big)& \frac{K(\bC(\{A\},H,\phi,n))}{K(\bC(Z^{c}\cap \{A\},H,\phi,n))}\ar[l]\ar[r]&\frac{K(\bC_{(n)}(\{A\} )}{K(\bC_{(n)}(Z^{c}\cap \{A\} ))}\\ 
K\big(\frac{C(Z\cap \{A\},H,\phi,n)}{C(Z\cap \{Z^{c}\}\cap\{A\},H,\phi,n)}\big)\ar[d] \ar[u] &\frac{K(\bC(Z\cap \{A\},H(Z),\phi_{|Z},n))}{K(\bC(Z\cap \{Z^{c}\}\cap\{A\},H(Z),\phi_{|Z},n))}\ar[u]\ar[l]\ar[d] \ar[r]&\frac{K(\bC_{(n)}(Z\cap \{A\}) )}{K(\bC_{(n)}(Z\cap \{Z^{c}\}\cap\{A\} ))}\ar[u]\ar[d] \\ 
K\big(\frac{C(Z^{\prime}\cap \{A^{\prime}\},H^{\prime},\phi^{\prime},n)}{C(Z^{\prime}\cap \{Z^{\prime,c}\}\cap\{A^{\prime}\},H^{\prime},\phi^{\prime},n)}\big)\ar[d] &\frac{K(\bC(Z^{\prime}\cap \{A^{\prime}\},H^{\prime},\phi^{\prime},n))}{K(\bC(Z^{\prime}\cap \{Z^{\prime,c}\}\cap\{A^{\prime}\},H^{\prime}(Z^{\prime}),\phi^{\prime}_{|Z^{\prime}},n))}\ar[l]\ar[d]\ar[r]&\frac{K(\bC_{(n)}(Z^{\prime}\cap \{A^{\prime}\} ))}{K(\bC_{(n)}(Z^{\prime}\cap \{Z^{\prime,c}\}\cap\{A^{\prime}\} ))}\ar[d]\\ 
K\big(\frac{C(\{A^{\prime}\},H^{\prime},\phi^{\prime},n)}{C(\{Z^{\prime,c}\}\cap\{A^{\prime}\},H^{\prime},\phi^{\prime},n )}\big)&\frac{K(\bC(\{A^{\prime}\},H^{\prime},\phi^{\prime},n))}{K(\bC(Z^{\prime,c}\cap \{A^{\prime}\},H^{\prime},\phi^{\prime},n))}\ar[l]\ar[r]&\frac{K(\bC_{(n)}(\{A^{\prime}\} ))}{K(\bC_{(n)}(Z^{\prime,c}\cap \{A^{\prime}\} ))}
}}}\notag\\
\mbox{}
\end{align}
 In the middle and right column we write cofibres of morphisms between $K$-theory spectra as quotients.
  The horizontal maps are induced from versions of the morphisms appearing in   \eqref{rgejknerkggergegergre2ge}. 
  They are equivalences by the proof of Theorem \ref{fwefiwjfeiooi234234324434en}. 
  The vertical morphisms in the upper and lower rows are  equivalences and the middle vertical morphisms are induced by the morphism $i$.
Let us explain how we get the left horizontal morphisms.
They all arise by the following general principle.
We consider a commutative diagram of graded $C^{*}$-algebras
$$\xymatrix{&C\ar[r]\ar[d]&D\ar[d]&&\\0\ar[r]&A\ar[r]&B\ar[r]& {B/A} \ar[r]&0} $$
Then we get  the commutative diagram of spectra
$$\xymatrix{K(C)\ar[d]\ar[r]&K(D)\ar[d]\ar[r]&\frac{K(D)}{K(C)}\ar@{..>}[d]\ar[r]&\Sigma K(C)\ar[d]\\K(A)\ar[r]&K(B)\ar[r]&K(B/A)\ar[r]&\Sigma K(A)}$$
where the filler of the left square and the fact that the lower sequence is a cofibre sequence yields the dotted arrow.

   By the definition of the excision morphism we have a commutative diagram
\begin{equation}\label{kopkopphrtopkpokrth}
\xymatrix{\frac{K(\bC_{(n)}(\{A\} ))}{K(\bC_{(n)}(Z^{c}\cap \{A\} ))}\ar[r]^(0.4){\simeq}\ar[d]&\KX_{(n)}^{G}(\{A\} ,\{Z^{c}\}\cap \{A\})\ar[d]^{e} \\ \frac{K(\bC_{(n)}(\{A^{\prime}\} ))}{K(\bC_{(n)}(Z^{\prime,c}\cap \{A^{\prime}\} ))}\ar[r]^(0.4){\simeq}&\KX_{(n)}^{G}(\{A^{\prime}\} ,\{Z^{\prime,c}\}\cap \{A^{\prime}\}) }
\end{equation}
where the left vertical {morphism} is the composition of the spectrum  morphisms in the right column
 of the diagram  \eqref{ferfpokp43rgr} from top to down.
 
 By  \eqref{knferlknlefnlrefref34r1} 
  the image of the relative index class of $\Dirac_{M}$ in the $(1,2)$-entry (we are counting from the left to the right, and from top to down) of \eqref{ferfpokp43rgr} maps to image of the relative index class of $\Dirac_{M^{\prime}}$ in the $(1,3)$-entry. 
 The class
 $\overline{\Ind(\Dirac_{M}, \on(A))}$ is obtained from the class in the $(1,2)$-entry  by going up and then right to the $(3,1)$-entry. Similarly, the class $\overline{\Ind(\Dirac_{M^{\prime}}, \on(A^{\prime}))}$ is obtained from the class in the $(1,3)$-entry by going down and then right to the  $(3,3)$-entry. 
   
 The equality asserted in the theorem now follows from the commutativity of the diagrams \eqref{ferfpokp43rgr}  and  \eqref{kopkopphrtopkpokrth}.
\end{proof}

\section{Suspension}

Let $G$ be a group and $M$ be a complete Riemannian manifold    with a proper   action of $G$ by isometries. We consider an invariant   Dirac operator $\Dirac$ on
$M$ of degree $n$. It acts on sections of a graded equivariant Dirac bundle $E\to M$  of right $\Cl^{n}$-modules. 

We now consider the Riemannian manifold $\tilde M:=\R\times M$ with the product metric $dt^{2}+g$.  Here $t$ is the coordinate of $\R$ and $g$ denotes the metric on $M$. The Riemannian manifold $\tilde M$ is complete and has an induced proper action of $G$ by isometries.

The pull-back $\tilde E^{\prime}\to \tilde M$ of the bundle $E\to M$ with the induced metric and connection is again equivariant.   We form the graded bundle $\tilde E:=\tilde E^{\prime}\otimes \Cl^{1}$. In view of the   identification $\Cl^{n+1}\cong \Cl^{n}\otimes \Cl^{1}$
it has a right action of the Clifford algebra  $\Cl^{n+1}$. The factor $\Cl^{n}$ acts on $\tilde E^{\prime}$, and the $\Cl^{1}$-factor acts on the $\Cl^{1}$-factor of $\tilde E$ by right-multiplication. 

The Clifford multiplication $TM\otimes E\to E$ extends to a Clifford multiplication $T\tilde M\otimes \tilde E\to \tilde E$, such that
$\partial_{t}$ acts by left-multiplication by the generator of $\Cl^{1}$ on the $\Cl^{1}$-factor.
In this way $\tilde E\to \tilde M$ becomes an invariant Dirac bundle of right $\Cl^{n+1}$-modules. We let
$\tilde \Dirac$ denote the associated  Dirac operator of degree $n+1$.
 
We assume that $\Dirac$ is uniformly positive outside of an invariant subset $A$ of $M$ {which is a support}. Then the operator $\tilde \Dirac$ will be uniformly positive outside of the subset $\R\times A$ of $\R\times M$. {Note that $\IR \times A$ is again a support.}
For every member $A^{\prime}$ of the big family $\{A\}$ the product 
  $\R\times A^{\prime}$ has an invariant coarsely excisive decomposition
$((-\infty,0]\times A^{\prime},[0,\infty)\times A^{\prime})$. It gives rise to a Mayer--Vietoris sequence.  
Since the entries of the decomposition with their induced structures are flasque
the boundary map of the  Mayer-Vietoris sequence is an  equivalence. Using the naturality of the Mayer--Vietoris sequences and  taking the colimit over the big family we get the equivalence of spectra
$$\partial:\KX_{(n+1)}^{G}(\R\times \{A\})\stackrel{\simeq}{\to} \Sigma \KX_{(n+1)}^{G}(\{A\})\ .$$
It induces an isomorphism of equivariant coarse $K$-homology groups
$$\delta:\pi_{0}\KX^{G}_{(n+1)}(\R\times \{A\})\to \pi_{-1}\KX^{G}_{(n+1)}(\{A\})\simeq \pi_{0} K\cX^{G}_{(n)}(\{A\})$$
such that $$\xymatrix{\pi_{0}\KX^{G}_{(n+1)}(\R\times \{A\})\ar[r]^-{\delta}\ar[d]^{\cong}&\pi_{0} K\cX^{G}_{(n)}(\{A\})\ar[d]^{\cong}\\ 
K\cX_{n+1}^{G}(\R\times \{A\})
\ar[r]^{\partial}&K\cX_{n}^{G}(  \{A\}}$$
commutes, where the vertical isomorphisms are given by Proposition \ref{kgowpegrfwfwfref}.

\begin{theorem}[Suspension]\label{fiweofwefwefewf}
We have the equality
$$\delta (\Ind(\tilde \Dirac,\on(\R\times A)))= \Ind(\Dirac,\on(A))\ .$$
\end{theorem}
\begin{proof}
The interesting analytic part of  the proof is Zeidler's \cite[Thm.~5.5]{MR3551834} showing the analogue of the assertion for the index classes in the $K$-theory of   Roe algebras associated to the situation. The rest is a tedious tour through various identifications made in order to interpret the index classes as equivariant coarse $K$-homology classes.

%

We choose an equivariant  ample $\tilde M$-controlled Hilbert space
$(\tilde H,\tilde \phi)$ and an equivariant ample $M$-controlled Hilbert space $(H,\phi)$.

Similarly as in the second part of the proof of Theorem \ref{roigeorgergergerg} we have a commutative diagram
\begin{equation}\label{giojigo345ggegergergergergre}
\xymatrix{K (C(\R\times \{A\},\tilde H,\tilde \phi))\ar[r]^{\simeq}\ar[d]&\KX^{G}(\R\times \{A \})\ar[d]\ar@/^3cm/[dddd]^{\partial}\\ 
\frac{K (C(\R\times \{A\},\tilde H,\tilde \phi))}{K (C([0,\infty)\times \{A\},\tilde H, \tilde \phi ))} \ar[r]^{\simeq}&\frac{\KX^{G}(\R\times \{A\})}{\KX^{G}([0,\infty)\times \{A\} )} \\ 
\frac{K (C([0,\infty)\times \{A\},\tilde H,\tilde \phi))}{K (C(\{\{0\}\}\times \{A\},\tilde H, \tilde \phi ))} \ar[u]^{\simeq}\ar[r]^{\simeq}\ar[d]^{\delta^{MV}} & \frac{\KX^{G}([0,\infty)\times \{A\} )}{\KX^{G}(\{\{0\}\}\times \{A\} )}\ar[d]^{\delta^{MV}}\ar[u]^{\simeq}\\\Sigma K (C(\{\{0\}\}\times \{A\},\tilde H, \tilde \phi ))\ar[r]^{\simeq}&\Sigma \KX^{G}(\{\{0\}\}\times \{A\} )\\
\Sigma K(C(\{A\},H,\phi))\ar[u]^{\simeq}\ar[r]^{\simeq}&\Sigma \KX^{G}(\{A\})\ar[u]^{\simeq}
}
\end{equation}
The horizontal maps are induced by   versions of \eqref{rgejknerkggergegergrege1}.
The lower left vertical morphism uses a controlled isometric embedding
$i_{*}(H,\phi)\to (\tilde H,\tilde \phi)$, where $i:M\cong \{0\}\times M \to \R\times M$
  is the embedding.  The two lower vertical morphisms are equivalences because they are induced by
a colimit over the  coarse equivalences $A^{\prime}\cong \{0\}\times A^{\prime} \to [0,n]\times A^{\prime}$ over $n$ in $\nat$ and  the members $A^{\prime}$ of the big family $\{A\}$ and $n$ in $\nat$. The filler of the bottom square comes from the second part of Theorem \ref{fwefiwjfeiooi234234324434e}.
  
  The morphism $\delta^{MV}$ in the left column is the boundary map in a fibre sequence of $K$-theory spectra induced from a long exact sequence of $C^{*}$-algebras.
 As the horizontal maps are eventually  induced from  zig-zags of morphisms between $C^{*}$-algebras fitting into respective sequences, and   the boundary map in the pair sequence for $\KX^{G}$
 comes from a boundary morphism associated to an exact sequence of $C^{*}$-algebras at the other end of the zig-zag  (the latter exact sequence is induced by applying the exact functor $A$ \cite[Prop. 8.9.2]{crosscat} to the exact sequence of $C^{*}$-categories determined by the ideal inclusion as in \cite[Lem. 8.83]{buen} or \cite[Lem. 6.16]{bu}), we get the filler of the square involving the boundary maps $\delta^{MV}$.

The choice of   controlled unitary embeddings $$(H_{0} ,\phi_{0} ,n )\to (\hat H,\phi,n)\ , \quad  
(\tilde H_{0},\tilde \phi_{0},n+1)\to (\widehat{\tilde  H},\widehat{\tilde \phi},n+1)$$
with ample complements  induces the   horizontal maps
 in the following commutative diagram:  
\begin{equation}\label{giojigo345ggegergergergergre1}
 \xymatrix{K_{0} (C(\R\times \{A\},\tilde H_{0},\tilde \phi_{0},n+1))\ar[r]\ar[d]\ar@/^-4cm/@{..>}[ddddd]^{d}&K_{0} (C(\R\times \{A\},\widehat{\tilde H},\widehat{\tilde \phi},n+1) \ar[d]\\
K_{0}\Big(\frac{ C(\R\times \{A\},\tilde H_{0},\tilde \phi_{0},n+1)}{C([0,\infty)\times \{A\},\tilde H_{0}, \tilde \phi_{0},n+1 )}\Big) \ar[r]&K_{0}\Big(\frac{ C(\R\times \{A\},\widehat{\tilde H},\widehat{\tilde \phi},n+1)}{C([0,\infty)\times \{A\},\widehat{\tilde H}, \widehat{\tilde \phi} ,n+1)}\Big)  \\
K_{0}\Big(\frac{   C([0,\infty)\times \{A\},\tilde H_{0},\tilde \phi_{0},n+1) }{   C(\{\{0\}\}\times \{A\},\tilde H_{0}, \tilde \phi_{0} ,n+1) }\Big) \ar[r] \ar[d]^{\delta^{MV}}\ar[u]^{\cong}&K_{0}\Big(\frac{   C([0,\infty)\times \{A\},\widehat{\tilde H},\widehat{\tilde \phi},n+1) }{   C(\{\{0\}\}\times \{A\},\widehat{\tilde H}, \widehat{\tilde \phi} ,n+1) }\Big)  \ar[d]^{\delta^{MV}} \ar[u]^{\cong}\\
K_{-1} (C(\{\{0\}\}\times \{A\}, \tilde H_{0},\ \tilde \phi_{0},n+1 ))\ar[r]&K_{-1} (C(\{\{0\}\}\times \{A\},\widehat{\tilde H},\widehat{ \tilde \phi},n+1 ))  \\
K_{-1}(C(\{A\}, H_{0}\otimes \Cl^{1}, \phi_{0}\otimes \id,n+1))\ar[r]\ar[u]^{\cong}&K_{-1}(C(\{A\},\hat H \otimes \Cl^{1},\hat \phi\otimes \id ,n+1))\ar[u]^{\cong}\\K_{0}(C(\{A\},H_{0},\phi_{0},n)\ar[u]^{\cong }\ar[r]& K_{0}(C(\{A\},\hat H,\hat \phi,n))\ar[u]^{\cong}
}\end{equation}

 The map induced in $K$-theory induced by the left column of \eqref{giojigo345ggegergergergergre} fits with the map induced by the right column in \eqref{giojigo345ggegergergergergre1} up to the isomorphisms of the kind
 $$K_{\ell+n}(C(\{A\},H,\phi))\cong K_{\ell}(C(\{A\},\hat H,\hat \phi,n))\cong K_{\ell-1}(C(\{A\},\hat H\otimes \Cl^{1},\hat \phi\otimes\id,n+1))$$ discussed at the end of Section \ref{goijrgoiregregreg}.

 Note that $i(\tilde \Dirac,\on(\R\times A))$ is an element in the upper left corner of \eqref{giojigo345ggegergergergergre1}. Its image in the upper right corner of
 \eqref{giojigo345ggegergergergergre} is $\Ind(\tilde \Dirac,\on(\R\times A))$ considered as a class in $K\cX_{n}^{G}(\R\times \{A\})$.
The class $i(\Dirac,\on(A))$ is a class in the lower left corner of  \eqref{giojigo345ggegergergergergre1}. Its image in the lower right corner of
 \eqref{giojigo345ggegergergergergre} is $\Ind( \Dirac, \on( A))$ considered as a class in $K\cX_{n}^{G}(\{A\})$.


Therefore, in order to show Theorem~\ref{fiweofwefwefewf}, by a diagram chase, it suffices to show that the dotted arrow $d$   satisfies $$d(i(\tilde \Dirac,\on(\R\times A)))=i(\Dirac,\on(A))\ .$$
 This equality follows  from the last assertion of Zeidler \cite[Thm.~5.5]{MR3551834} by applying the evaluation map from the localization algebra to the Roe algebra. Zeidler's proof only uses Proposition \ref{gropgergergerge} in order to reduce to the special case $M=\R$ (with trivial action). He assumes free $G$-actions, but this is not relevant for this part of the argument.
\end{proof}

\bibliographystyle{alpha}
\bibliography{indexona}

\end{document}